\documentclass{amsart}

\newtheorem{theorem}{Theorem}[section]
\newtheorem{lemma}[theorem]{Lemma}

\newtheorem{corollary}[theorem]{Corollary}
\newtheorem{proposition}[theorem]{Proposition}
\newtheorem{claim}[theorem]{Claim}
\usepackage{graphicx}

\theoremstyle{definition}

\theoremstyle{remark}

\numberwithin{equation}{section}

\begin{document}

\title[Dehn fillings on large hyperbolic $3$-manifolds]
  {Toroidal Dehn fillings on large hyperbolic $3$-manifolds}

\author{Masakazu Teragaito}
\address{Department of Mathematics and Mathematics Education, Hiroshima University,
1-1-1 Kagamiyama, Higashi-hiroshima, Japan 739-8524}
\email{teragai@hiroshima-u.ac.jp}
\thanks{
Partially supported by Japan Society for the Promotion of Science,
Grant-in-Aid for Scientific Research (C), 16540071.
}%



\dedicatory{Dedicated to Professor Takao Matumoto on his sixtieth birthday}


\begin{abstract}
We show that if a hyperbolic $3$-manifold $M$ with a single torus boundary admits two Dehn fillings
at distance $5$, each of which contains an essential torus, then $M$ is a rational homology solid torus, which is not
large in the sense of Wu.
Moreover, one of the surgered manifold contains an essential torus which meets the core of the attached solid torus
minimally in at most two points.
This completes the determination of best possible upper bounds for the distance between two exceptional Dehn fillings yielding
essential small surfaces in all ten cases for large hyperbolic $3$-manifolds.
\end{abstract}

\maketitle

\section{Introduction}\label{sec:1}

Let $M$ be a compact orientable $3$-manifold with a torus boundary component $T_0$.
A \textit{slope\/} on $T_0$ is the isotopy class of an essential unoriented simple closed curve on $T_0$.
For two slopes $\gamma_1,\gamma_2$, the \textit{distance\/} $\Delta(\gamma_1,\gamma_2)$ between them is their minimal geometric
intersection number.
For a slope $\gamma$ on $T_0$,
the manifold obtained by \textit{$\gamma$-Dehn filling\/} on $M$ is $M(\gamma)=M\cup V$, where
$V$ is a solid torus glued to $M$ along $T_0$ in such a way that $\gamma$ bounds a meridian disk in $V$.
If $M$ is hyperbolic, in the sense that $M$ with its boundary tori removed admits a complete hyperbolic structure
with totally geodesic boundary, then a slope $\gamma$, or the filling, is said to be \textit{exceptional\/} if
$M(\gamma)$ is not hyperbolic.
In particular, if $M(\gamma)$ contains an essential torus, then $\gamma$, or the filling, is said to be \textit{toroidal}.
We are interested in obtaining the upper bounds for the distance between exceptional slopes, and
focus on toroidal slopes in this paper.

Following Wu \cite{W2},
let us say that $M$ is \textit{large\/} if $H_2(M,\partial M-T_0)\ne 0$.
Note that $M$ is not large if and only if $M$ is a $\mathbb{Q}$-homology solid torus or a $\mathbb{Q}$-homology cobordism between two tori. 
Hence, $M$ is large if $\partial M$ contains at least three tori or a component of genus at least two.
Wu \cite{W2} showed that the upper bound for the distance between exceptional fillings can be often improved when
we restrict ourselves to large hyperbolic $3$-manifolds.
For example, the distance between toroidal fillings on a hyperbolic $3$-manifold with a torus boundary component is at most $8$, but
it is at most $5$ for large hyperbolic $3$-manifolds \cite{Go2}, because
the only hyperbolic manifolds with a pair of toroidal fillings at distance greater than $5$ are 
obtained by Dehn fillings on the Whitehead link exterior, so that they are all $\mathbb{Q}$-homology solid tori.
In \cite[Question 4.2]{Go4}, Gordon asks if there is a large hyperbolic manifold with toroidal fillings at distance $5$.
In this direction, \cite[Theorem 3.1]{BGZ} shows that if $\partial M$ is a single torus and the first betti number $\beta_1(M)\ge 3$ then
the distance between two toroidal fillings is at most $4$.
As stated in \cite[Remark 3.15]{BGZ}, their argument also works for $M$ whose boundary consists of at least $4$ tori.
In \cite{GT}, we showed that if $M$ admits two toroidal filling at distance $5$, then $\partial M$ consists of either a single torus or two tori.
Furthermore, Lee \cite{L} proved that the Whitehead sister link (the $(-2,3,8)$-pretzel link) exterior, which is not large, is the only 
manifold whose boundary consists of two tori and which admits two toroidal fillings at distance $5$.

In this paper, we analyze the case where $M$ has a single torus boundary and show the following.

\begin{theorem}\label{thm:main}
Let $M$ be a hyperbolic $3$-manifold whose boundary is a single torus $T_0$.
If there are two toroidal slopes $\alpha$ and $\beta$ on $T_0$ with $\Delta(\alpha,\beta)=5$,
then $M$ is a $\mathbb{Q}$-homology solid torus.
Moreover, either $M(\alpha)$ or $M(\beta)$ contains an essential torus which meets the core of the attached solid torus
minimally in at most two points.
\end{theorem}

There are infinitely many examples as in Theorem \ref{thm:main}.
In particular, the exteriors of Eudave-Mu\~{n}oz knots $k(2,-1,n,0)$ $(n\ne 1)$ \cite{EM} give
all knot exteriors in $S^3$ that satisfy the condition of Theorem \ref{thm:main} \cite{T2}.
In these examples, each of surgered toroidal manifolds contains an essential torus which meets the core
of the attached solid torus in two points.
%

As a corollary, we can answer to Gordon's question \cite[Question 4.2]{Go4}.

\begin{corollary}\label{cor:main}
Let $M$ be a large hyperbolic $3$-manifold with a torus boundary component $T_0$.
If $M$ admits two toroidal filling $\alpha$ and $\beta$ on $T_0$, then $\Delta(\alpha,\beta)\le 4$.
\end{corollary}

As mentioned in \cite{Go4}, this upper bound is sharp.  For example, the Whitehead link exterior is large, and
admits toroidal slopes $0$ and $4$ on one boundary torus.
Corollary \ref{cor:main} completes the determination of best possible upper bounds for the distance between two exceptional Dehn fillings yielding
essential small surfaces in all ten cases for large hyperbolic $3$-manifolds.  These are shown in Table \ref{table:large}, where
$S$, $D$, $A$ and $T$ indicate that the manifold $M(\alpha)$ or $M(\beta)$ contains an essential sphere, disk, annulus or torus, respectively.
For these bounds, refer to \cite{FMP,Go4}.

\begin{table}[ht]
\renewcommand\arraystretch{1.5}
\noindent\[
\begin{array}{|c||c|c|c|c|}
\hline
\Delta & S & D & A & T \\
\hline\hline
S & 0 & 0 & 1 & 1 \\
\hline
D & & 1 & 2 & 1 \\
\hline
A & & & 4 & 4 \\
\hline T & & & & 4 \\
\hline
\end{array}
\]
\caption{Upper bounds on $\Delta(\alpha,\beta)$ for large hyperbolic $3$-manifolds}\label{table:large}
\end{table}

Combining with known facts \cite{Go2,L}, we have:

\begin{corollary}\label{cor:hit}
If a hyperbolic $3$-manifold $M$ with a torus boundary component $T_0$ admits two toroidal slopes on $T_0$ at distance
greater than four, then either surgered manifold contains an essential torus which meets the core of
the attached solid torus minimally in at most two points.
\end{corollary}

The proof of Theorem \ref{thm:main} goes as follows.
Assume that $M$ admits two toroidal fillings at distance $5$.
First, we deal with the case where both fillings yield no Klein bottle.
As shown in \cite[Proposition 2.3]{GT}, we may assume that at least one essential torus is separating.
By following the arguments of \cite{GT}, we will see that this torus meets the core of the attached solid torus only twice. 
The argument is divided into three cases, according to how many times the other essential torus meets the core of the attached solid
torus.
Secondly, we consider the case where either surgered manifold contains a Klein bottle.
The proof of Theorem \ref{thm:main} will be completed in Section \ref{sec:proofs},
where the proofs of Corollaries \ref{cor:main} and \ref{cor:hit} are also given.
The proof of Theorem \ref{thm:main} implicitly gives the collection of $3$-manifolds that includes
all hyperbolic $3$-manifolds with a single torus boundary
such that there are two toroidal fillings at distance $5$.

We use integer coefficients for homology groups.

\section{Preliminaries}

Let $M$ be a hyperbolic $3$-manifold whose boundary is a torus $T_0$.
Suppose that $M$ admits two toroidal slopes $\alpha$ and $\beta$ on $T_0$ with $\Delta(\alpha,\beta)=5$.
Then $M(\alpha)$ and $M(\beta)$ are irreducible by \cite{O,W}.

Let $\widehat{S}$ be an essential torus in $M(\alpha)$.
By \cite[Proposition 2.3]{GT}, we may assume that $\widehat{S}$ is separating.
We may assume that $\widehat{S}$ meets the attached solid torus $V_\alpha$ in $s$ disjoint meridian disks $u_1,u_2,\dots,u_s$,
numbered successively along $V_\alpha$ and that $s$ is minimal over all choices of $\widehat{S}$.
Let $S=\widehat{S}\cap M$.
By the minimality of $s$, $S$ is incompressible and boundary-incompressible in $M$.
Similarly, we choose an essential torus $\widehat{T}$ in $M(\beta)$ which meets the attached solid torus $V_\beta$ in $t$ disjoint
meridian disks $v_1,v_2,\dots,v_t$, numbered successively along $V_\beta$, where $t$ is minimal as above.
Thus we have another incompressible and boundary-incompressible torus $T=\widehat{T}\cap M$.
We may assume that $S$ and $T$ intersect transversely.
Then $S\cap T$ consists of arcs and circles.  Since both surfaces are incompressible, we can assume that $S\cap T$ contains no circle component
bounding a disk in $S$ or $T$.
Moreover, we can assume that $\partial u_i$ meets $\partial v_j$ in $5$ points for any pair of $i$ and $j$.
Orient all boundary components $\partial u_i$ of $S$ coherently on $T_0$.
Similarly, orient all $\partial v_j$ of $T$ coherently on $T_0$.
We can choose an oriented meridian-longitude pair $m$ and $l$ on $T_0$ so that $[\partial u_i]=[m]$ and $[\partial v_j]=d[m]+5[l]$ for some $d$ in $H_1(T_0)$.
Furthermore, we can assume that $d=1$ or $2$ by reversing the orientations of all $\partial v_j$ and $l$ if necessary \cite{Go2}.
This number $d$ is called the \textit{jumping number\/} of $\alpha$ and $\beta$.

\begin{lemma}\label{lem:jumping}
Let $a_1,a_2,a_3,a_4,a_5$ be the points of $\partial u_i\cap\partial v_j$, numbered so that they appear successively on $\partial u_i$ along its orientation.
Then these points appear in the order of $a_d,a_{2d},a_{3d},a_{4d},a_{5d}$ on $\partial v_j$ along its orientation.
In particular, 
if $d=1$, then
two points of $\partial u_i\cap\partial v_j$ are successive on $\partial u_i$ if and only if they are successive on $\partial v_j$,
and if $d=2$, then
two points of $\partial u_i\cap\partial v_j$ are successive on $\partial u_i$ if and only if they are not successive on $\partial v_j$.
\end{lemma}

\begin{proof}
This immediately follows from the definition of the jumping number.
See \cite[Lemma 2.10]{GW}.
\end{proof}

\begin{figure}[tb]
\includegraphics*[scale=0.6]{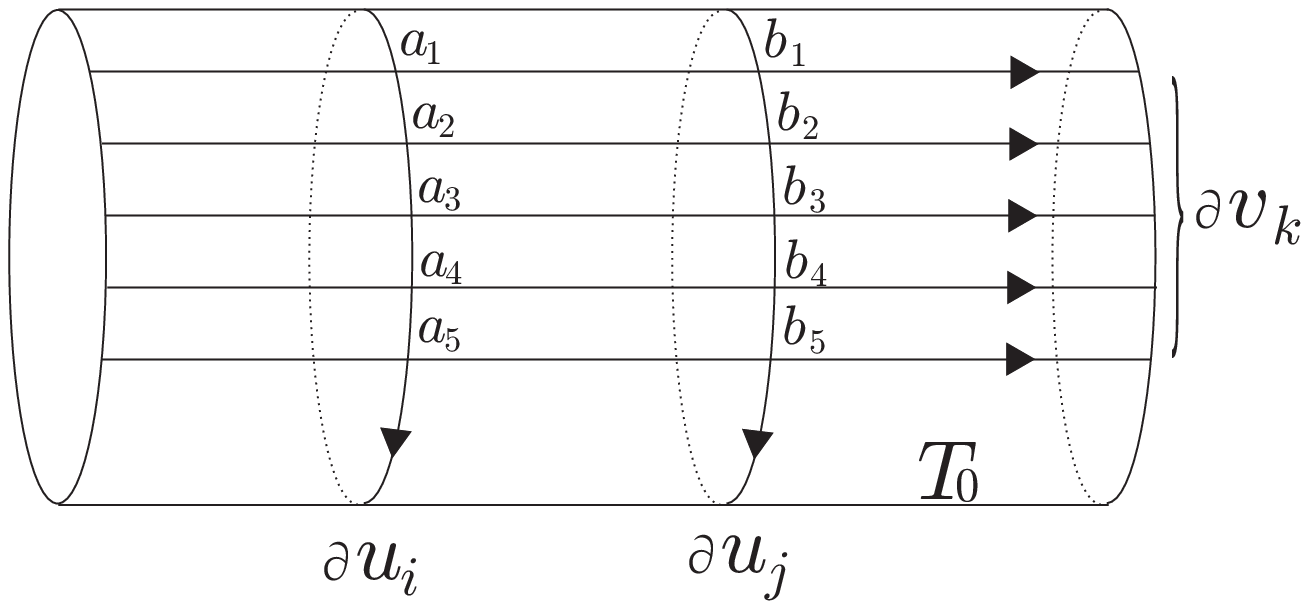}
\caption{}\label{fig:jump}
\end{figure}

This simple lemma is important, and will be used repeatedly in the paper.
For example, let $\partial u_i\cap \partial v_k=\{a_1,a_2,a_3,a_4,a_5\}$
and $\partial u_j\cap \partial v_k=\{b_1,b_2,b_3,b_4,b_5\}$ as shown in Figure \ref{fig:jump}.
If the jumping number is one, then
these points appear in the order $a_1$,$b_1$,$a_2$,$b_2$,$a_3$,$b_3$,$a_4$,$b_4$,$a_5$,$b_5$ on $\partial v_k$ along its orientation.
If the jumping number is two, then 
they appear in the order $a_1$,$b_1$,$a_3$,$b_3$,$a_5$,$b_5$,$a_2$,$b_2$,$a_4$,$b_4$ on $\partial v_k$ along its orientation.

Let $G_S$ be the graph on $\widehat{S}$ consisting of the $u_i$ as (fat) vertices and the arcs of $S\cap T$ as edges.
The orientation on $\partial u_i$ induces that of $u_i$.
Thus each vertex is assigned a sign.
Define $G_T$ on $\widehat{T}$ similarly.
Two graphs on a surface are considered to be equivalent if there is a homeomorphism of the surface carrying one graph to the other.
Note that both graphs have no trivial loops, since $S$ and $T$ are boundary-incompressible.

For an edge $e$ of $G_S$ incident to $u_i$,
the endpoint of $e$ is labelled $j$ if it is in $\partial u_i\cap \partial v_j$.
Similarly, label the endpoints of each edge of $G_T$.
Thus the labels $1,2,\dots,t$ (resp.\ $1,2,\dots,s$) appear in order around each vertex of $G_S$ (resp.\ $G_T$) repeated $5$ times.
Since $S$ is orientable, we can distinguish the signs of $u_i$'s, according as 
the labels $1,2,\dots,t$ appear counterclockwise or clockwise, if $t>2$.
The situation for $G_T$ is similar.
Each vertex $u_i$ of $G_S$ has degree $5t$, and each $v_j$ of $G_T$ has degree $5s$.
If an edge $e$ has labels $j_1,j_2$ at its endpoints, then $e$ is called a \textit{$\{j_1,j_2\}$-edge}.

Let $G=G_S$ or $G_T$.
An edge of $G$ is a \textit{positive\/} edge if it connects vertices of the same sign.
Otherwise it is a \textit{negative\/} edge.
Possibly, a positive edge is a loop.

A cycle in $G$ consisting of positive edges is a \textit{Scharlemann cycle\/} if it bounds a disk face of $G$ and all edges
in the cycle are $\{i,i+1\}$-edges for some label $i$.
The number of edges in a Scharlemann cycle is called the \textit{length\/} of the Scharlemann cycle, and
the set $\{i,i+1\}$ is called its \textit{label pair}.
A Scharlemann cycle of length two is called an \textit{$S$-cycle} for short.

\begin{lemma}\label{lem:common}
\begin{itemize}
\item[(1)] (The parity rule) An edge is positive in a graph if and only if it is negative in the other graph. 
\item[(2)] There are no two edges which are parallel in both graphs.
\item[(3)] The edges of a Scharlemann cycle of $G_S$ \textup{(}resp.\ $G_T$\textup{)} cannot lie in a disk on $\widehat{T}$ \textup{(}resp.\ $\widehat{S}$\textup{)}.
\end{itemize}
\end{lemma}

\begin{proof}
(1) See \cite[p.279]{CGLS}.
(2) is \cite[Lemma 2.1]{Go2}.
For (3), see \cite[Lemma 2.2(5)]{GW}.
\end{proof}

Let $e_1,e_2,\dots,e_t$ be $t$ mutually parallel negative edges in $G_S$ numbered successively, each connecting vertex $u_i$ to $u_j$.
Suppose that $e_k$ has label $k$ at $u_i$ for $1\le k\le t$.
Then this family defines a permutation $\sigma$ of the set $\{1,2,\dots,t\}$ such that
$e_k$ has label $\sigma(k)$ at $u_j$.
In fact, $\sigma(k)\equiv k+h \pmod{t}$ for some $h$.
We call $\sigma$ the \textit{associated permutation\/} to the family. 
It is well-defined up to inversion.

\begin{lemma}\label{lem:GS}
$G_S$ satisfies the following.
\begin{itemize}
\item[(1)] If $G_S$ contains a Scharlemann cycle, then $\widehat{T}$ is separating.
\item[(2)] Let $t\ge 3$.  Any family of mutually parallel positive edges contains at most $t/2+1$ edges.
If it contains more than $t/2$ edges, then it contains an $S$-cycle.
\item[(3)] Let $t\ge 3$. If a family of mutually parallel negative edges contains more than $t$ edges,
then all the vertices of $G_T$ have the same sign, and the associated permutation to this family has a single orbit.
\item[(4)] If $t\ge 4$, then any family of mutually parallel edges contains at most $2t$ edges.
\end{itemize}
\end{lemma}

\begin{proof}
For (1), see \cite[Lemma 2.2(4)]{GW}.
(2) is essentially the same as \cite[Lemma 2.8]{T}.
(3) is \cite[Lemma 2.3(1)]{GW}.
(4) is \cite[Corollary 5.5]{Go2}.
\end{proof}

The next lemma will be used repeatedly throughout the paper.

\begin{lemma}\label{lem:homology}
If $H_1(M(\alpha))$ \textup{(}or $H_1(M(\beta))$\textup{)} is finite, then $M$ is a $\mathbb{Q}$-homology solid torus.
\end{lemma}

\begin{proof}
If $H_1(M(\alpha))$ is finite, then $\beta_1(M)=1$ by Poincar\'{e} duality and the Mayer-Vietoris sequence.
\end{proof}

\section{No Klein bottle}\label{sec:noklein}

Until the end of Section \ref{sec:t3}, we assume that 
neither $M(\alpha)$ nor $M(\beta)$ contains a Klein bottle.

\begin{lemma}\label{lem:s2}
$s=2$.
\end{lemma}

\begin{proof}
Assume $s\ge 4$.
By \cite[Proposition 4.1]{GT}, $t\ge 2$.
Moreover, Section 5 of \cite{GT} shows $t\ne 2$.
Thus $t\ge 3$.
Then Section 3 of \cite{GT} eliminates the case $t\ge 3$.
\end{proof}

\begin{figure}[tb]
\includegraphics*[scale=0.7]{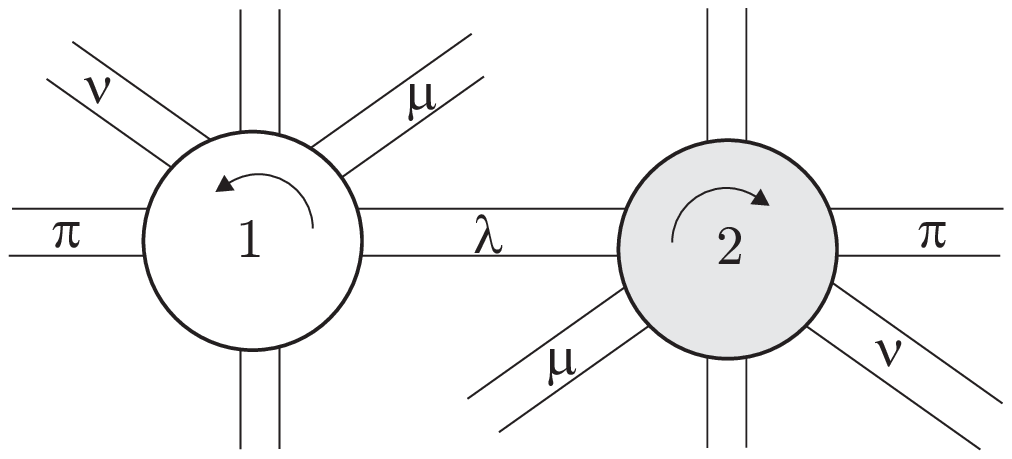}
\caption{}\label{fig:edgelabel}
\end{figure}

\begin{figure}[tb]
\includegraphics*[scale=0.5]{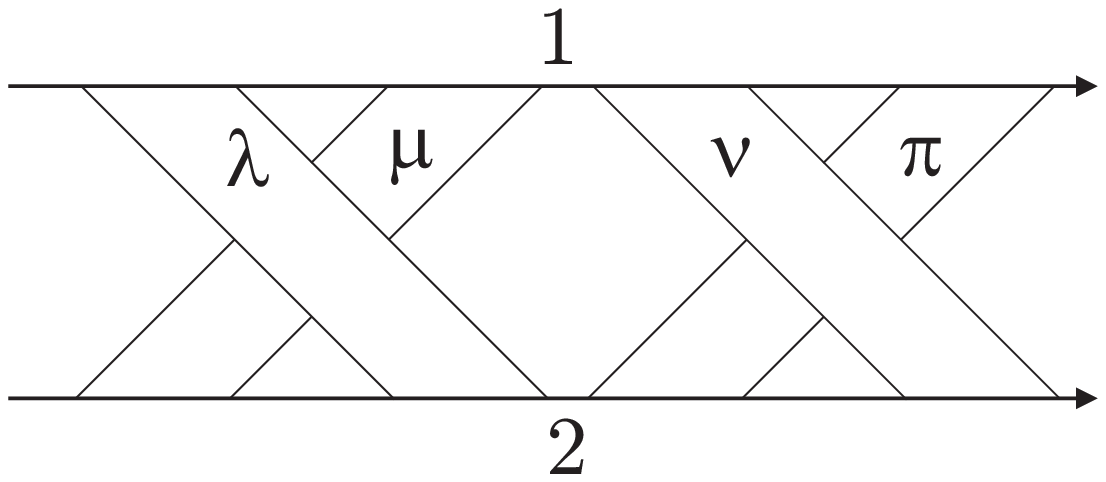}
\caption{}\label{fig:edgeclass}
\end{figure}

Thus the non-loop edges of $G_S$ are divided into at most four classes $\lambda$, $\mu$, $\nu$ and $\pi$, called \textit{edge classes\/}, of
mutually parallel edges.  See Figure \ref{fig:edgelabel}.
The orientation of each vertex is indicated by an arrow inside the vertex throughout the paper.
Also, $G_S$ are described more schematically as in Figure \ref{fig:edgeclass}, where loops are disregarded.
An edge $e$ of $G_T$ is labelled by the class of the corresponding edge of $G_S$, which is referred to as the \textit{edge class label\/} of $e$.

Two of edge classes are said to be \textit{adjacent\/} if the endpoints at $u_1$ (hence at $u_2$ as well) of those edge classes
are successive.
Otherwise, they are \textit{non-adjacent}.
Thus if all edge classes are not empty, then the pair $\{\lambda,\nu\}$ and the pair $\{\mu,\pi\}$ are non-adjacent.

Let $\varepsilon,\delta$ be two elements of $\{\lambda,\mu,\nu,\pi\}$.
A face $f$ of $G_T$ is called an \textit{$(\varepsilon,\delta)$-face\/} if any edge on $\partial f$
has an edge class label $\varepsilon$ or $\delta$.
Since such $f$ is bounded by a Scharlemann cycle, each label appears at least once by Lemma \ref{lem:common}(3).
An $(\varepsilon,\delta)$-face $f$ is \textit{good\/}
if either label $\varepsilon$ or $\delta$ does not appear consecutively along $\partial f$.
We remark that any bigons and $3$-gons bounded by positive edges in $G_T$ are good (see \cite[Lemma 3.7]{GL4} for $3$-gons).


We collect some results here which will be used in the following sections.

Let $\alpha_1$, $\alpha_2$, $\alpha_3$, $\alpha_4$ be the number of edges in the edge classes $\lambda,\mu,\nu,\pi$, respectively. 
Also, let $\alpha_0$ be the number of loops at each vertex of $G_S$.
(Clearly, the two vertices are incident to the same number of loops.)
Then $G_S$ is determined by a quintuple $(\alpha_0,\alpha_1,\alpha_2,\alpha_3,\alpha_4)$.
We say then $G_S\cong G(\alpha_0,\alpha_1,\alpha_2,\alpha_3,\alpha_4)$.
If $\alpha_0=0$, then we abbreviate it to $G(\alpha_1,\alpha_2,\alpha_3,\alpha_4)$.
Note that 
\begin{align*}
G(\alpha_0,\alpha_1,\alpha_2,\alpha_3,\alpha_4) & \cong G(\alpha_0,\alpha_2,\alpha_1,\alpha_4,\alpha_3)\cong G(\alpha_0,\alpha_3,\alpha_4,\alpha_1,\alpha_2) \\
 & \cong G(\alpha_0,\alpha_4,\alpha_3,\alpha_2,\alpha_1),
\end{align*}
and that
$G(\alpha_1,\alpha_2,\alpha_3,\alpha_4)\cong G(\alpha_2,\alpha_3,\alpha_4,\alpha_1)$ (see \cite{Go2}).

Each edge class corresponds to either loops in $G_T$, or non-loop edges in $G_T$.
Define $\varepsilon_i$ to be $0$ or $1$ according as the edge class with $\alpha_i$ edges is of the first or second kind.

\begin{lemma}\label{lem:epsilon}
If $\varepsilon_i=0$ then $\alpha_i\le 2$, and if $\varepsilon_i=1$ then $\alpha_i\le 4$.
Moreover, $\alpha_i+\varepsilon_i\equiv \alpha_j+\varepsilon_j \pmod{2}$ for $i,j\in \{1,2,3,4\}$.
\end{lemma}

\begin{proof}
This is \cite[Lemma 5.3]{Go2}.
\end{proof}

Since $\widehat{S}$ is separating, it divides $M(\alpha)$ into a \textit{black side\/} $\mathcal{B}$ and a \textit{white side\/} $\mathcal{W}$.
Also, the faces of $G_T$ are divided into \textit{black\/} and \textit{white\/} faces, according as they lie in $\mathcal{B}$ or $\mathcal{W}$.

\begin{lemma}\label{lem:atmost3}
Suppose that all the vertices of $G_T$ have the same sign.  Then $G_T$ satisfies the following.
\begin{itemize}
\item[(1)] Two bigons with the same color have the same pair of edge class labels.
\item[(2)] At most three edges can be mutually parallel.
\item[(3)] A black bigon and a white bigon cannot have the same pair of edge class labels.  Also,
if their pairs are disjoint, then each pair consists of non-adjacent edge class labels.
\item[(4)] If a bigon has non-adjacent edge class labels, then there is no $3$-gon with the same color as the bigon.
If a bigon has adjacent edge class labels, $\{\lambda,\mu\}$ say, then
any $3$-gon with the same color as the bigon has edge class labels $\{\nu,\pi\}$. 
\item[(5)] If there are a black bigon and a white bigon with disjoint edge class labels, then there is no $3$-gon.
\item[(6)] There cannot be good $(\varepsilon,\delta)$-faces for both colors.
\end{itemize}
\end{lemma}

\begin{proof}
(1) This is \cite[Lemma 5.2]{GL3}.  (Recall that $M(\alpha)$ does not contain a Klein bottle.)

(2) Assume that $G_T$ contains mutually parallel four edges.
Then there are two bigons with the same color among these edges.
But (1) and Lemma \ref{lem:common}(2) give a contradiction.

(3) Let $f$ be a black bigon and $g$ be a white bigon.
Let $H=V_\alpha\cap \mathcal{B}$.
Then shrinking $H$ to its core in $H\cup f$ gives a M\"{o}bius band $B$ in $\mathcal{B}$ whose boundary lies on $\widehat{S}$.
The same procedure gives another M\"{o}bius band $B'$ in $\mathcal{W}$.
If $f$ and $g$ have the same edge class label pair, or if each pair consists of adjacent edge class labels and
the two pairs are disjoint, then $\partial B$ is isotopic to $\partial B'$ on $\widehat{S}$, and hence
$M(\alpha)$ contains a Klein bottle.

(4) We may assume that a black bigon $f$ has the edge class label pair $\{\lambda,\nu\}$.
If there is a black $3$-gon $g$, then its edge class labels are $\{\lambda,\nu\}$ by \cite[Lemma 3.7]{GL5},
Then there is an essential annulus $A$ in $\widehat{S}$ which contains all edges of $f$ and $g$.
Let $F$ be a torus obtained from $A$ by capping $\partial A$ off with two disks.
On $F$, $\partial f$ and $\partial g$ give disjoint essential loops which represent different homology classes, a contradiction.
If $f$ has $\{\lambda,\mu\}$, then $g$ has $\{\lambda,\mu\}$ or $\{\nu,\pi\}$ by \cite[Lemma 3.7]{GL5}.
However the former is impossible as above.

(5) follows from (3) and (4).

(6) Let $f$ be a good black $(\varepsilon,\delta)$-face and $g$ a good white $(\varepsilon,\delta)$-face.
Then there is an essential annulus $A$ in $\widehat{S}$ which contains all edges of $f$ and $g$.
Let $X=N(A\cup H\cup f)\subset \mathcal{B}$, where $H=V_\alpha\cap \mathcal{B}$.
Then $X$ is a solid torus, in which the core of $A$ is homotopic to at least twice the core of $X$ by \cite[Lemma 4.4]{GL5}.
Let us write $\mathcal{B}=X\cup X'$.
Then $X'$ is also a solid torus (see the proof of \cite[Theorem 4.1]{GL5}).
Let $Y=N(A\cup H'\cup g)\subset \mathcal{W}$, where $H'=V_\alpha\cap \mathcal{W}$, and let us write $\mathcal{W}=Y\cup Y'$.
Similarly, $Y$ and $Y'$ are solid tori.
Thus $X\cup Y$ is a Seifert fibered manifold over the disk with two exceptional
fibers, and a regular fiber is given by the core of $\partial A$.
Also, $X'\cup Y'$ also admits such a Seifert fibration.
Thus $\partial (X\cup Y)$ gives an essential torus in $M(\alpha)$ which is disjoint from $V_\alpha$.
This contradicts the minimality of $s$.
\end{proof}

\section{The case where $t\le 2$}\label{sec:t2}

First, we assume $t=1$.

\begin{proposition}\label{pro:t1}
If $t=1$, then $M$ is a $\mathbb{Q}$-homology solid torus.
\end{proposition}

\begin{figure}[tb]
\includegraphics*[scale=0.6]{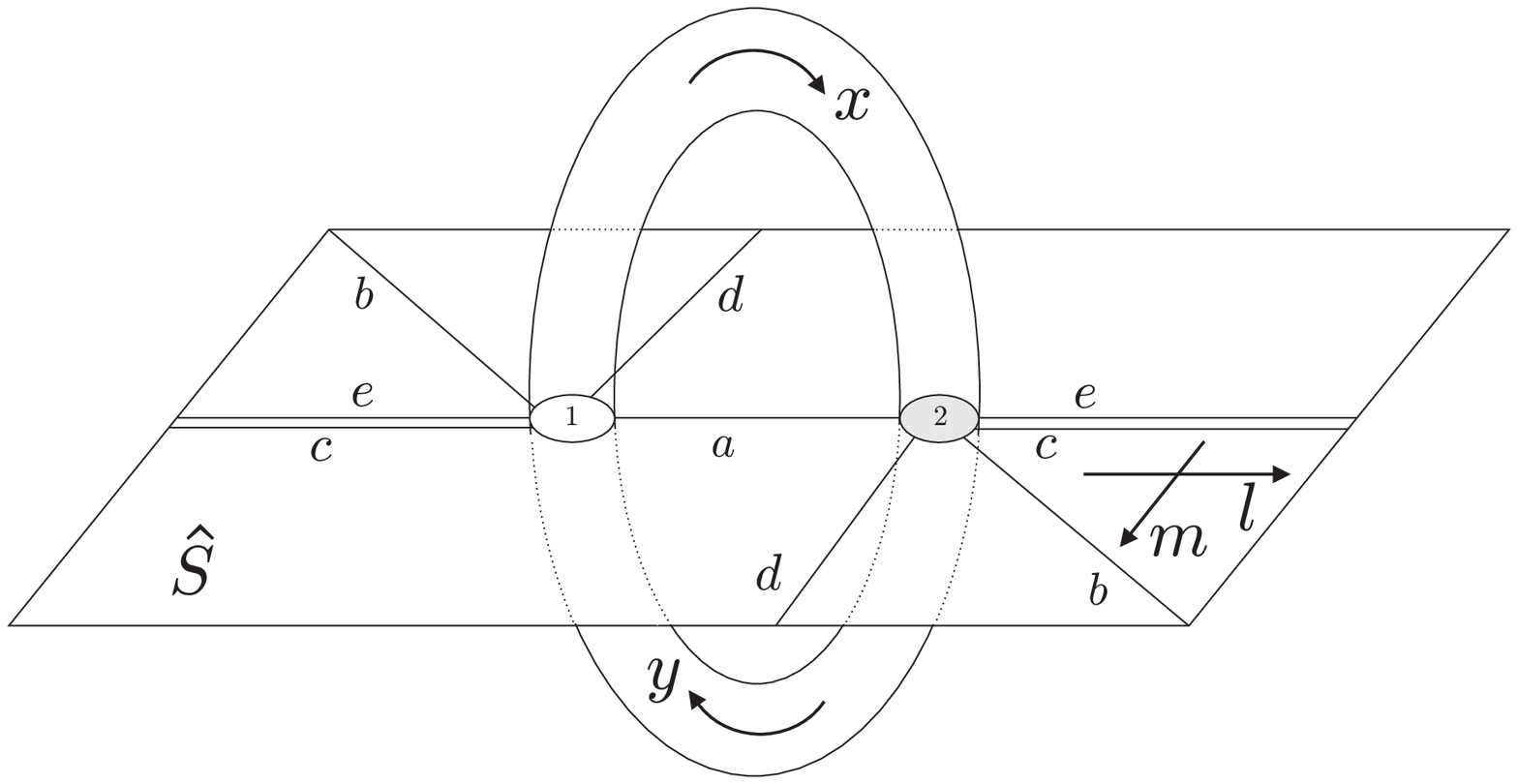}
\caption{}\label{fig:t1s2-hom}
\end{figure}

\begin{proof}
There are only two possible graph pairs as shown in Figures \ref{fig:t1s2} and \ref{fig:t1s2-2} by \cite[Section 4]{GT}.
The jumping number is one and two, respectively.

\begin{figure}[tb]
\includegraphics*[scale=0.6]{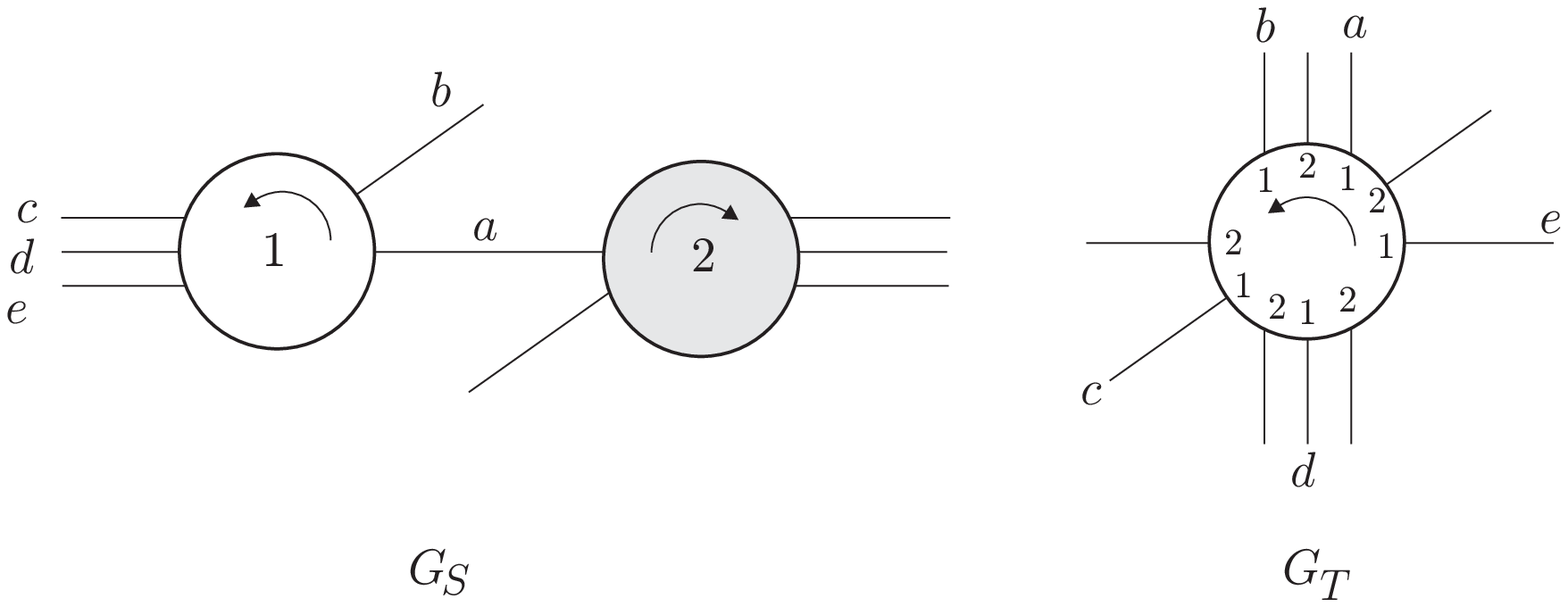}
\caption{}\label{fig:t1s2}
\end{figure}

\begin{figure}[tb]
\includegraphics*[scale=0.6]{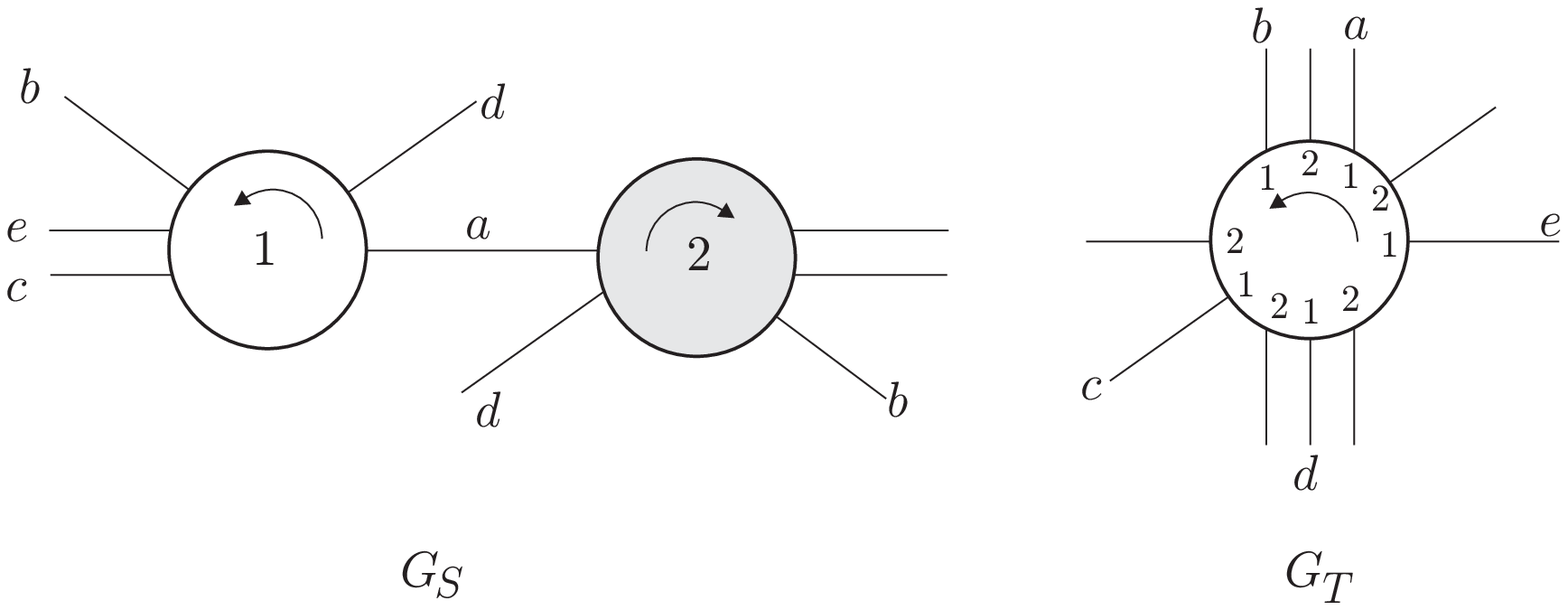}
\caption{}\label{fig:t1s2-2}
\end{figure}

We show that Figure \ref{fig:t1s2} is impossible.
Let $f_1$ be the bigon bounded by $\{a,d\}$
and $f_2$ be the $3$-gon bounded by $\{b,c,e\}$ in $G_T$.
They lie on the same side of $\widehat{S}$.
We see that $f_1$ has edge class labels $\{\lambda,\nu\}$, and $f_2$ has edge class labels $\{\mu,\nu\}$.
These two sets are distinct and not disjoint.
This is impossible by Lemma \ref{lem:atmost3}(4).

We calculate $H_1(M(\alpha))$ when the jumping number is two.
Let $f_1$ and $g_1$ be the bigons bounded by $\{a,d\}$ and $\{d,b\}$, respectively.
Also, let $f_2$ and $g_2$ be the $3$-gons bounded by $\{b,c,e\}$ and $\{a,c,e\}$, respectively.
Then $f_1$ and $f_2$ lie on the same side of $\widehat{S}$.
We may assume that this side is $\mathcal{B}$.
Let $H=V_\alpha\cap \mathcal{B}$.
Then $\mathcal{B}=N(\widehat{S}\cup H\cup f_1\cup f_2)\cup (\text{a $3$-ball})$.
For, $\partial f_1$ is non-separating on the genus two surface $F$ obtained from $\widehat{S}$ by tubing along $H$,
and $\partial f_2$ is non-separating on the torus obtained from $F$ by compressing along $f_1$.
Since $M(\alpha)$ is irreducible, 
the $2$-sphere obtained by compression along $f_2$
bounds a ball in $\mathcal{B}$.
The situation in $\mathcal{W}$ is similar (use $g_1$, $g_2$ instead of $f_1$, $f_2$).

Take a generator $\ell,m,x,y$ of $H_1(\widehat{S}\cup V_\alpha)$ as in Figure \ref{fig:t1s2-hom}, where
$x$ is represented by the union of the core of the upper half part $H$ of $V_\alpha$ and the edge $a$,
and $y$ is similar.
Then we have $[\partial f_1]= 2x+m$, $[\partial f_2]= 3x+3\ell+m$,
$[\partial g_1]= 2y-2m-\ell$ and $[\partial g_2]=3y-2\ell$.
Thus $H_1(M(\alpha))$ has a presentation 
\[\langle \ell,m,x,y\ |\ 2x+m=0, 3x+3\ell+m=0, 2y-2m-\ell=0, 3y-2\ell=0\rangle.
\]
It is easy to show that $H_1(M(\alpha))=\mathbb{Z}_{35}$.
Thus $M$ is a $\mathbb{Q}$-homology solid torus by Lemma \ref{lem:homology}. 
\end{proof}

Now we assume $t=2$.
Then $G_T$ is determined by a quintuple $(\beta_0,\beta_1,\beta_2,\beta_3,\beta_4)$ as $G_S$,
and then we say $G_T\cong G(\beta_0,\beta_1,\beta_2,\beta_3,\beta_4)$.

\begin{lemma}\label{lem:t2opposite}
If the two vertices of $G_T$ have distinct signs, then $M$ is a $\mathbb{Q}$-homology solid torus.
\end{lemma}

\begin{proof}
Assume that the two vertices of $G_T$ have distinct signs.
Then there is only one possible graph pair as shown in Figure \ref{fig:s2t2}
by \cite[Lemmas 7.3 and 7.4]{GT}.
Note that the jumping number is two.
We remark that 
$G_S$ can be either graph of Figure \ref{fig:s2t2} and
$\widehat{T}$ is separating as well $\widehat{S}$ under this situation by Lemma \ref{lem:GS}(1).

\begin{figure}[tb]
\includegraphics*[scale=0.5]{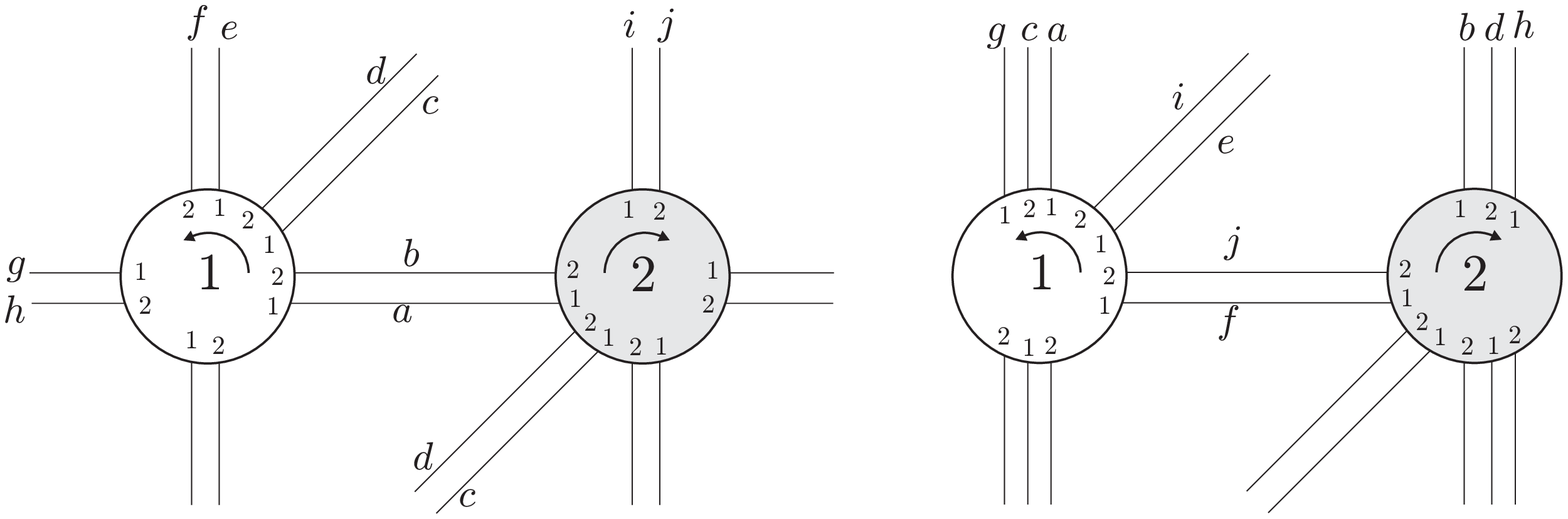}
\caption{}\label{fig:s2t2}
\end{figure}

Without loss of generality, we can assume that $G_S$ is the first one of Figure \ref{fig:s2t2}.
As in the proof of Proposition \ref{pro:t1}, we calculate $H_1(M(\alpha))$.
Let $f_1$ be the bigon bounded by $\{a,c\}$ in $G_T$. We may assume that $f$ is black.
Let $H=V_\alpha\cap\mathcal{B}$.
There is an essential annulus $A$ on $\widehat{S}$ which contains the edges $a,c$.
Let $X=N(A\cup H\cup f_1)$.
Then $X$ is a solid torus, where the core of $A$ runs twice along the core of $X$ \cite[Lemma 3.7]{GL2}.
Let us write $\mathcal{B}=X\cup X'$.
Then $X'$ is also a solid torus by the minimality of $\widehat{S}$, and
moreover the core of the annulus $\widehat{S}-A$ runs at least twice along the core of $X'$. (See the proof of \cite[Theorem 3.2]{GL2}.)
For the other side $\mathcal{W}$ of $\widehat{S}$,
let $g_1$ be the bigon bounded by $\{c,g\}$ and $g_2$ the $3$-gon bounded by $\{a,f,i\}$ in $G_T$.
Then $\mathcal{W}=N(\widehat{S}\cup H'\cup g_1\cup g_2)\cup (\text{a $3$-ball})$, where $H'=V_\alpha\cap \mathcal{W}$, as in the proof of Proposition \ref{pro:t1}.

By using the generators $x,y,\ell,m$ of $H_1(\widehat{S}\cup V_\alpha)$ as shown in Figure \ref{fig:t1s2-hom},
$[\partial f_1]=2x+m$, $[\partial g_1]=2y-m-\ell$, and $[\partial g_2]=2m+y$.
Notice that $\ell+x$ and $m$ generate $H_1(\partial X')$.
Thus if a meridian of $X'$ represents $p(\ell+x)+qm$, then $|p|\ge 2$.
Hence $H_1(M(\alpha))$ has a presentation
\[\langle \ell,m,x,y\ |\ 2x+m=0, p(\ell+x)+qm=0, 2y-m-\ell=0, 2m+y=0\rangle,
\]
which is equivalent to $\langle x\ | \ (11p-2q)x=0\rangle$.
Hence if $p/q\ne 2/11$, then $M$ is a $\mathbb{Q}$-homology solid torus by Lemma \ref{lem:homology}.

Suppose $p/q=2/11$. 
There is a M\"{o}bius band $B$ properly embedded in $X$ whose boundary is homotopic to the core of $A$
as in the proof of Lemma \ref{lem:atmost3}(3).
Since $p=\pm 2$, the curve $m$, which is parallel to the core of $A$, runs twice along the core of $X'$.
This implies that it bounds M\"{o}bius band $B'$ in $X'$.
Then $M(\alpha)$ contains a Klein bottle obtained from $B\cup B'$, contradicting our assumption.
\end{proof}


Thus we consider the case where the two vertices of $G_T$ have the same sign.
Then $G_T$ contains only positive edges, so $G_S$ contains only negative edges by the parity rule.
Hence the edges of $G_S$ are divided into at most four edge classes, and so $G_S\cong G(\alpha_1,\alpha_2,\alpha_3,\alpha_4)$.
Since $\sum_{i=1}^{4}\alpha_i=10$, the possibilities for $(\alpha_1,\alpha_2,\alpha_3,\alpha_4)$ allowed by Lemma \ref{lem:epsilon} are
$(4,4,2,0)$, $(4,4,1,1)$, $(4,1,4,1)$, $(4,2,2,2)$, $(3,3,3,1)$, $(3,3,2,2)$ and $(3,2,3,2)$, up to equivalence.

\begin{lemma}\label{lem:4420}
$(4,4,2,0)$ is impossible.
\end{lemma}

\begin{proof}
By Lemma \ref{lem:epsilon}, $\varepsilon_1=\varepsilon_2=\varepsilon_3=1$.
Thus $G_T$ contains no loops, so the edges of $G_T$ are also divided into at most four edge classes.
The edges of the class $\lambda$ of $G_S$ belong to mutually distinct edge classes by Lemma \ref{lem:common}(2).
This also holds for the edges of $\mu$ and $\nu$.
Hence $G_T\cong G(0,3,3,2,2)$ or $G(0,3,2,3,2)$.
However there is no correct labeling for either of the configurations.
\end{proof}

\begin{lemma}
$(4,4,1,1)$ and $(4,1,4,1)$ are impossible.
\end{lemma}

\begin{proof}
Consider the case $(4,4,1,1)$.
By Lemma \ref{lem:epsilon}, $\varepsilon_1=\varepsilon_2=1$ and $\varepsilon_3=\varepsilon_4=0$.
Then $G_T\cong G(1,2,2,2,2)$.
Thus $G_T$ contains a black bigon and a white bigon with the same pair of edge class labels $\{\lambda,\mu\}$,
contradicting Lemma \ref{lem:atmost3}(3).
$(4,1,4,1)$ is ruled out in the same way.
\end{proof}

\begin{lemma}\label{lem:4222}
$(4,2,2,2)$ is impossible.
\end{lemma}

\begin{proof}
All $\varepsilon_i$'s are $1$ by Lemma \ref{lem:epsilon}.
By Lemma \ref{lem:atmost3}(2), $G_T\cong G(0,3,3,3,1)$, $G(0,3,3,2,2)$ or $G(0,3,2,3,2)$.

If $G_T\cong G(0,3,3,3,1)$, then there are three bigons with the same color.
By Lemma \ref{lem:atmost3}(1), they have the same pair of edge class labels.
This means that $\alpha_i\ge 3$ for at least two $i$'s, a contradiction.

For the remaining two cases, there is no correct labeling in $G_T$.
\end{proof}

\begin{lemma}
If $(\alpha_1,\alpha_2,\alpha_3,\alpha_4)=(3,3,3,1)$, then $M$ is a $\mathbb{Q}$-homology solid torus.
\end{lemma}

\begin{proof}
All $\varepsilon_i$'s are $1$ again.
As in the proof of the previous lemma, $G(0,3,3,3,1)$ is the only possibility of $G_T$. 
Then we can determine the correspondence between the edges of $G_S$ and $G_T$ as in Figure \ref{fig:3331},
where the jumping number is two.

\begin{figure}[tb]
\includegraphics*[scale=0.6]{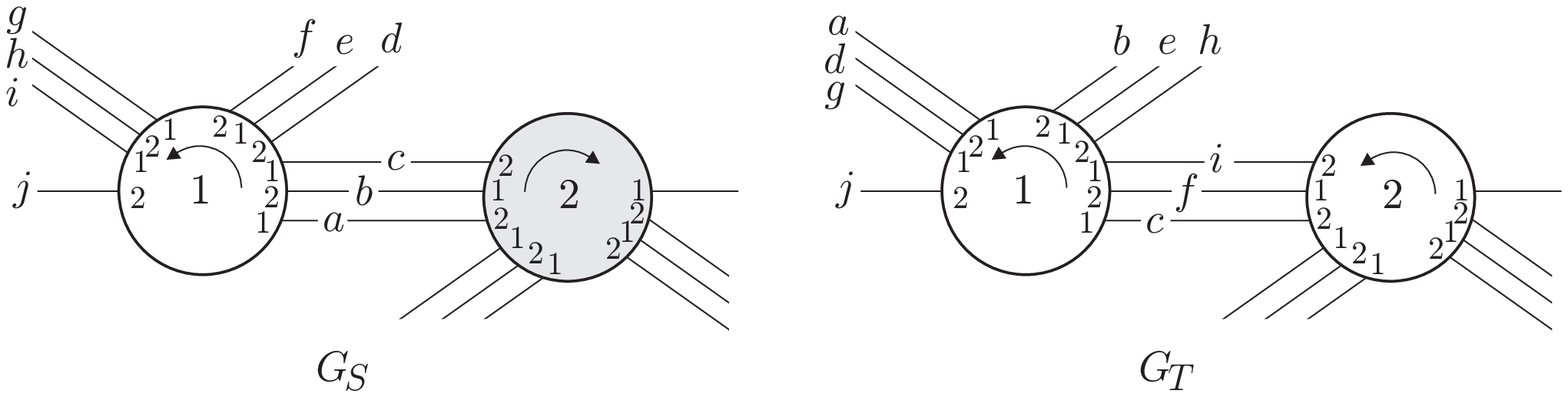}
\caption{}\label{fig:3331}
\end{figure}

We calculate $H_1(M(\alpha))$ as in the proof of Proposition \ref{pro:t1}.
Let $f_1$ be the bigon bounded by $\{c,f\}$ and $f_2$ the $4$-gon bounded by $\{g,h,i,j\}$ in $G_T$.
They lie on the same side $\mathcal{B}$, say, of $\widehat{S}$.
Then we see that $\mathcal{B}=N(\widehat{S}\cup H\cup f_1\cup f_2)\cup (\text{a $3$-ball})$, where $H=V_\alpha\cap \mathcal{B}$.
Also, let $g_1$ be the white bigon bounded by $\{f,i\}$ and $g_2$ the white $4$-gon bounded by $\{a,b,c,j\}$ in $G_T$.
Again, $\mathcal{W}=N(\widehat{S}\cup H'\cup g_1\cup g_2)\cup (\text{a $3$-ball})$, where $H'=V_\alpha\cap \mathcal{W}$.
Take the same generators $\ell,m,x,y$ of $H_1(\widehat{S}\cup V_\alpha)$ as in Figure \ref{fig:t1s2-hom}.
Then $[\partial f_1]=2x+m$, $[\partial f_2]=4x+4\ell+3m$, $[\partial g_1]=2y-\ell-2m$ and $[\partial g_2]=4y-\ell$.
Hence $H_1(M(\alpha))$ has a presentation 
\[
\langle \ell, m, x, y \ |\ 2x+m=0, 4x+4\ell+3m=0, 2y-\ell-2m=0, 4y-\ell=0\rangle.
\]
This shows $H_1(M(\alpha))=\mathbb{Z}_2\oplus\mathbb{Z}_{30}$.
Hence $M$ is a $\mathbb{Q}$-homology solid torus.
\end{proof}

\begin{lemma}\label{lem:t2last}
$(3,3,2,2)$ and $(3,2,3,2)$ are impossible.
\end{lemma}

\begin{proof}
Assume $G_S\cong G(3,3,2,2)$.
Then $\varepsilon_1=\varepsilon_2=1$ and $\varepsilon_3=\varepsilon_4=0$.
Thus each vertex of $G_T$ is incident to exactly two loops which are parallel.
By Lemma \ref{lem:atmost3}(2), $G_T\cong G(2,3,3,0,0)$, $G(2,3,1,2,0)$, $G(2,3,1,1,1)$, $G(2,2,2,2,0)$ or $G(2,2,2,1,1)$. 
For $G(2,3,3,0,0)$, $G(2,3,1,2,0)$ and $G(2,3,1,1,1)$, $G_T$ contains a black bigon and a white bigon with the same pair of edge class
labels $\{\lambda,\mu\}$, contradicting Lemma \ref{lem:atmost3}(3).
For $G(2,2,2,2,0)$ and $G(2,2,2,1,1)$, $G_T$ contains at least four bigons with the same color.
But their pairs of edge class labels are distinct, contradicting Lemma \ref{lem:atmost3}(1).

$(3,2,3,2)$ is ruled out in exactly the same way.
\end{proof}

\begin{proposition}\label{pro:t2}
If $t=2$, then $M$ is a $\mathbb{Q}$-homology solid torus.
\end{proposition}

\begin{proof}
If the two vertices of $G_T$ have distinct signs, then $M$ is a $\mathbb{Q}$-homology solid torus by Lemma \ref{lem:t2opposite}.
Otherwise, Lemmas \ref{lem:4420}-\ref{lem:t2last} shows that $G_S\cong G_T\cong G(0,3,3,3,1)$ and $M$ is a $\mathbb{Q}$-homology solid torus.
\end{proof}

\section{The case where $t\ge 3$}\label{sec:t3}

In this section, we consider the case where $t\ge 3$, which will be ruled out.
If $\widehat{T}$ is separating, then $t\ge 4$.
Then the argument of \cite[Section 5]{GT}, with exchanging the role of $S$ and $T$, eliminates this case.
Hence we may assume that $\widehat{T}$ is non-separating.

\begin{lemma}\label{lem:noloop}
$G_S$ contains no loops.
\end{lemma}

\begin{proof}
Assume $\alpha_0>0$.  Then $G_T$ contains a negative edge, and hence not all the vertices of $G_T$ have
the same sign.
Thus $\alpha_i\le t$ for $i=1,2,3,4$ by Lemma \ref{lem:GS}(3).
Since $\widehat{T}$ is non-separating,  $\alpha_0\le t/2$ by Lemma \ref{lem:GS}(1) and (2).
Hence $\sum_{i=1}^{4}\alpha_i \ge 4t$, and so $G_S\cong G(t/2,t,t,t,t)$.
But \cite[Lemma 5.2]{GT} (with an exchange of $G_S$ and $G_T$) eliminates this configuration.
\end{proof}

\begin{lemma}
All the vertices of $G_T$ have the same sign.
\end{lemma}

\begin{proof}
By the previous lemma, $\sum_{i=1}^{4}\alpha_i =5t$.
Hence some $\alpha_i>t$, and so all the vertices of $G_T$ have the same sign by Lemma \ref{lem:GS}(3).
\end{proof}

Thus every edge of $G_T$ is a positive $\{1,2\}$-edge, and every disk face of $G_T$ is a Scharlemann cycle with label pair $\{1,2\}$.
Let $D_2$ and $D_3$ be the number of bigons and $3$-gons of $G_T$, respectively.


\begin{lemma}\label{lem:D2}
$G_T$ has only disk faces. Moreover,
if $D$ is the number of disk faces of $G_T$,
then $D=4t$, $D_2\ge 2t$ and $2D_2+D_3\ge 6t$.
Moreover, if $D_2=2t$, then $G_T$ has only bigons and $3$-gons, and $D_3=2t$.
\end{lemma}

\begin{proof}
We may assume $\alpha_1>t$. By Lemma \ref{lem:GS}(3), the associated permutation $\sigma$ to the class $\lambda$ has a single orbit.
The first and $(t+1)$-th edges of $\lambda$ have the same label pair at its endpoints.
These two edges are not parallel in $G_T$ by Lemma \ref{lem:common}(2),
and hence cutting $\widehat{T}$ along the first $t+1$ edges of $\lambda$ gives a disk.
Hence $G_T$ has only disk faces.

Then Euler's formula $t-5t+D=0$ gives $D=4t$.
Also, $2D_2+3(D-D_2)\le 10t$ gives $D_2\ge 2t$.
Similarly, $2D_2+3D_3+4(D-D_2-D_3)\le 10t$ gives $2D_2+D_3\ge 6t$.

If $D_2=2t$, then $D_3\ge 2t$.
Thus $D\ge D_2+D_3\ge 4t$.
Since $D=4t$, we have $D=D_2+D_3$ and $D_3=2t$.
In particular, $G_T$ has only bigons and $3$-gons.
\end{proof}


The next is the key lemma in this section.

\begin{lemma}\label{lem:BW}
$G_T$ contains a black bigon and a white bigon.
\end{lemma}

\begin{proof}
By Lemma \ref{lem:D2}, $D_2\ge 2t$.
We divide the proof into two cases.

(1) $D_2=2t$.

Without loss of generality, we suppose that all bigons of $G_T$ are black.
Then all white faces are $3$-gons by Lemma \ref{lem:D2}.

If black bigons have non-adjacent edge class labels, then
all black faces are bigons by Lemma \ref{lem:atmost3}(4).
This means that any edge of $G_T$ belong to a bigon, and hence
there are only $2D_2=4t$ edges in $G_T$, a contradiction.
Hence we can assume that black bigons have the pair $\{\lambda,\mu\}$ of adjacent edge class labels.

We claim that there is a black $3$-gon.
For, if not, any $3$-gon is white.
Thus there are $2t$ white $3$-gons.
This implies that $G_T$ has at least $6t$ edge, a contradiction.
In fact, any black $3$-gon has edge class labels $\{\nu,\pi\}$ by Lemma \ref{lem:atmost3}(4).
Hence $\alpha_1=\alpha_2=2t$ and $\alpha_3+\alpha_4=t$.

Let $x$ and $y$ be the number of black $3$-gons with edge class labels $(\nu,\nu,\pi)$ and $(\nu,\pi,\pi)$, respectively.
Counting $\nu$- and $\pi$-edges, $(2x+y)+(x+2y)=t$, giving $x+y=t/3$.
Thus there are $D_3-t/3=5t/3$ white $3$-gons.

On the other hand, if a white $3$-gon has non-adjacent edge class labels,
then every white faces ($3$-gons) have such labels, and hence there are only two edge class labels in $G_T$, a contradiction.
Thus any white $3$-gon has adjacent edge class labels, so
it is either $\{\lambda,\pi\}$ or $\{\mu,\nu\}$ by Lemma \ref{lem:atmost3}(6).
Let $a,b,c,d$ be the number of white $3$-gons with edge class labels $(\lambda,\lambda,\pi)$, $(\lambda,\pi,\pi)$,
$(\mu,\mu,\nu)$ and $(\mu,\nu,\nu)$, respectively.
Counting $\nu$- and $\pi$-edges, $(c+2d)+(a+2b)=t$.
Since $a+b+c+d=5t/3$, we have $b+d=-2t/3 < 0$, a contradiction.

(2) $D_2>2t$.

We may assume that all bigons are black.
Let $\{\varepsilon,\delta\}$ be the pair of edge class labels of black bigons.
Then the edge class $\varepsilon$ (and $\delta$) of $G_S$ contains more than $2t$ edges.
Hence $t=3$ by Lemma \ref{lem:GS}(4), and so $D_2>6$.
Since $G_T$ has just $5t=15$ edges and any two bigons do not share an edge, $D_2\le 7$.
Therefore $D_2=7$.  This means that two of $\alpha_i$'s are at least $7$.
Since $\sum_{i=1}^{4}\alpha_i=15$, there are only two possibilities for $(\alpha_1,\alpha_2,\alpha_3,\alpha_4)$:
$(7,7,1,0)$ and $(7,8,0,0)$, up to equivalence.

Let $\nu(v_i,v_j)$ denote the number of mutually non-parallel edges in $G_T$ that join 
$v_i$ and $v_j$.

Suppose that $(\alpha_1,\alpha_2,\alpha_3,\alpha_4)=(7,7,1,0)$.
We may choose the labels around vertex $u_1$ as in Figure \ref{fig:7710}.
Then there are three possibilities for the labeling at $u_2$.

\begin{figure}[tb]
\includegraphics*[scale=0.6]{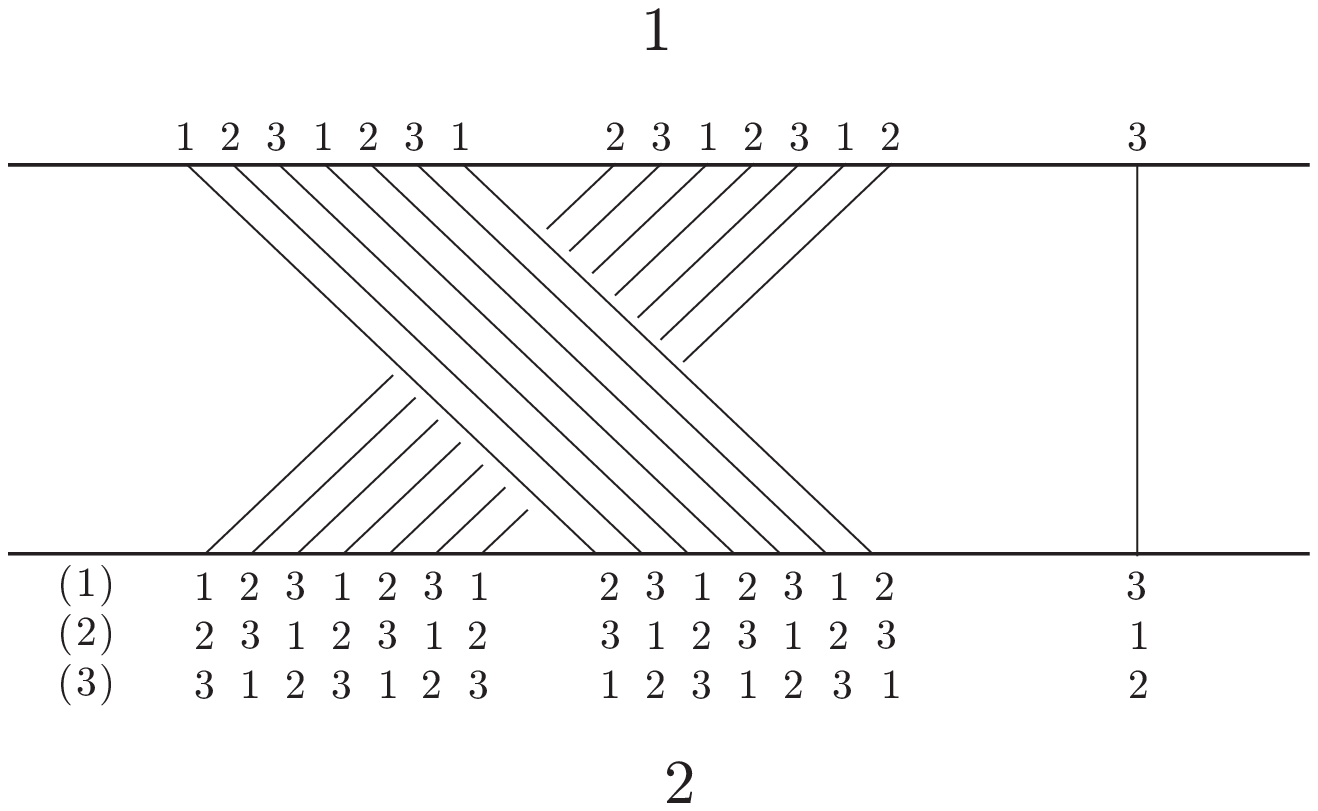}
\caption{}\label{fig:7710}
\end{figure}

For (1), $G_T$ contains a loop at $v_3$.
Then $\nu(v_1,v_2)\le 2$ by \cite[Lemma 7.4(i)]{Go2}.
However, the class $\lambda$ contains three $\{1,2\}$-edges.
Since they are not mutually parallel in $G_T$ by Lemma \ref{lem:common}(2),
$\nu(v_1,v_2)\ge 3$, a contradiction.
Similar arguments rule out (2) and (3).

Next, consider the case $(7,8,0,0)$ as shown in Figure \ref{fig:7800}.

\begin{figure}[tb]
\includegraphics*[scale=0.6]{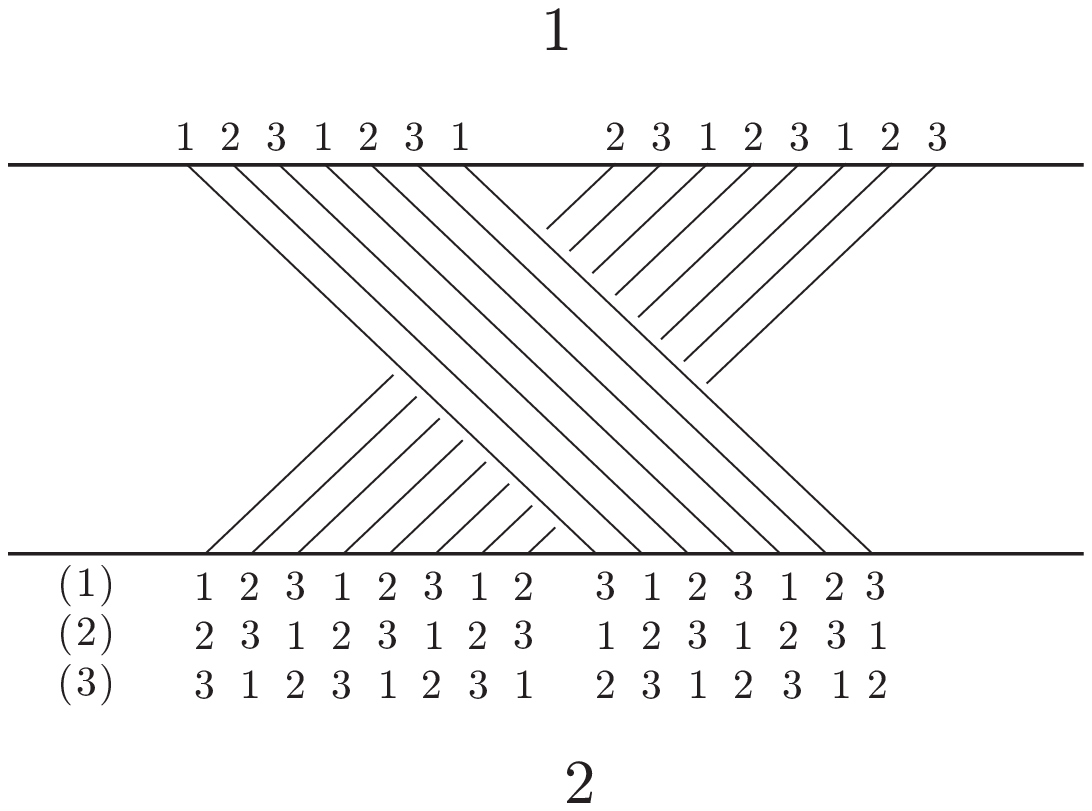}
\caption{}\label{fig:7800}
\end{figure}

For (1), $\nu(v_1,v_2)\ge 3$ by considering the edges of $\mu$.
Let $A,B,C$ be the three $\{1,3\}$-edges in $\lambda$, and let $D$ be one of $\{1,3\}$-edges in $\mu$.
Then we can see that  $D$ is not parallel in $G_T$ to either $A$, $B$ or $C$ by \cite[Lemma 2.5(i)]{Go2}.
(See also the proof of \cite[Lemma 7.8]{Go2}.)
Thus $\nu(v_3,v_1)\ge 4$, contradicting \cite[Lemma 7.4(ii)]{Go2}.
A similar argument rules out (3).
For (2), each vertex of $G_T$ is incident to a loop.
The class $\lambda$ contains three $\{1,1\}$-edges, and then these give
three mutually parallel loops at $v_1$.
This contradicts Lemma \ref{lem:common}(2).
(See also \cite[Lemma 7.6(i)]{Go2}.)
\end{proof}

\begin{lemma}\label{lem:bigon-common}
The pair of edge class labels of black bigons in $G_T$ has exactly one common label with the pair of edge class labels of white bigons.
In particular, at least one pair consists of adjacent edge class labels.
\end{lemma}

\begin{proof}
By Lemma \ref{lem:atmost3}(3), black bigons and white bigons cannot have the same pair of edge class labels.
Also, if black bigons and white bigons have disjoint pairs, then
both pairs consist of non-adjacent edge class labels.
Then there is no $3$-gon by Lemma \ref{lem:atmost3}(5).
Hence $D_2\ge 3t$ by Lemma \ref{lem:D2}.
Since any two bigons do not share an edge, this implies that $G_T$ contains at least $6t$ edges, a contradiction.
Thus black bigons and white bigons have a single common edge class label.
\end{proof}

\begin{proposition}\label{pro:tge3}
Case $t\ge 3$ is impossible.
\end{proposition}

\begin{proof}
By Lemma \ref{lem:bigon-common}, we may assume that
black bigons have the pair of adjacent edge class labels $\{\lambda,\mu\}$.
There are two cases.

(1) White bigons have the pair of adjacent edge class labels.

Then we can assume that they have the pair $\{\lambda,\pi\}$.
Thus $\lambda$ is the common label of black and white bigons.
Since $D_2\ge 2t$, $\alpha_1\ge 2t$, and $\alpha_2+\alpha_4\ge 2t$.
Hence $\alpha_3\le t$.

First, assume $D_2=2t$.  Then $D_3=2t$ by Lemma \ref{lem:D2}.
By Lemma \ref{lem:atmost3}(4), any black $3$-gon has edge class labels $\{\nu,\pi\}$ and any white $3$-gon
has edge class labels $\{\mu,\nu\}$.
Thus $\alpha_1=2t$, since there are only bigons and $3$-gons by Lemma \ref{lem:D2}.
Let $a$ and $b$ be the number of black $3$-gons and white $3$-gons, respectively.
Then $D_3=a+b$.
Notice that any black $3$-gon contains at least one edge with label $\nu$.
Since two black $3$-gons do not share an edge, $a\le \alpha_3\le t$.
Similarly, we have $b\le t$.
Thus $a=b=t$, and in fact, any black $3$-gon has edge class labels $(\nu,\pi,\pi)$ and any white $3$-gon
has edge class labels $(\mu,\mu,\nu)$.
Then any edge of $G_T$ with edge class label $\pi$ belongs to a black $3$-gon.
Hence $\alpha_4=2a=2t$.
Similarly, any edge with edge class label $\mu$ belongs to a white $3$-gon, and hence $\alpha_2=2b=2t$.
Then $\alpha_2+\alpha_4=4t$, a contradiction.

Next, assume $D_2>2t$.
Then $\alpha_1>2t$.  Hence $t=3$ by Lemma \ref{lem:GS}(4).
Thus $D_2>6$.
If $D_2\ge 8$, then $\alpha_1\ge 8$ and $\alpha_2+\alpha_4 \ge 8$.
This is impossible, since $\sum_{i=1}^{4}\alpha_i=15$.
Hence $D_2=7$, and so $\alpha_1\ge 7$, $\alpha_2+\alpha_4\ge 7$ and $\alpha_3\le 1$.
Again, any black $3$-gon has edge class labels $\{\nu,\pi\}$ and any white $3$-gon has edge class labels $\{\mu,\nu\}$.
Since any black (and white) $3$-gon contains at least one edge with edge class label $\nu$,
there is at least one $3$-gon for each color.
Hence $D_3\le 2$.
But $D_3\ge 18-2D_2=4$ by Lemma \ref{lem:D2}, a contradiction.

(2) White bigons have the pair of non-adjacent edge class labels.

We may assume that they have the pair $\{\lambda,\nu\}$.
Then $\alpha_1\ge 2t$, $\alpha_2+\alpha_3\ge 2t$ and hence $\alpha_4\le t$.

If $D_2=2t$, then $G_T$ contains only bigons and $3$-gons by Lemma \ref{lem:D2}.
But there is no white $3$-gon by Lemma \ref{lem:atmost3}(4).
Hence all white faces are bigons.
This means that any edge of $G_T$ belongs to a white bigon, and hence any edge has label $\lambda$ or $\nu$, a contradiction.
Thus $D_2>2t$, so $t=3$ as in (1).
Then $D_2=7$, and hence $\alpha_1\ge 7$, $\alpha_2+\alpha_3\ge 7$ and $\alpha_4\le 1$.
Since any $3$-gon is black and has edge class labels $\{\nu,\pi\}$, $D_3\le 1$.
This contradicts Lemma \ref{lem:D2} again.
\end{proof}

\section{Klein bottle case}\label{sec:klein}

From this section, we deal with the case where $M(\alpha)$ or $M(\beta)$ contains a Klein bottle.
Without loss of generality, we may assume that $M(\alpha)$ contains a Klein bottle.

\begin{lemma}\label{lem:beta-noklein}
If $M(\beta)$ also contains a Klein bottle, then $M=W(-4)$, where $W$ is the Whitehead link exterior. 
In particular, $M$ is a $\mathbb{Q}$-homology solid torus.
Moreover, both of $M(\alpha)$ and $M(\beta)$ contain Klein bottles meeting the core of the attached solid torus
once.
\end{lemma}

\begin{proof}
This is  \cite[Theorem 1.4]{L}.
\end{proof}

Let $\widehat{P}$ be a Klein bottle in $M(\alpha)$.
We may assume that $\widehat{P}\cap V_\alpha$ consists of
$p$ meridian disks $u_1,u_2,\dots,u_p$, numbered successively, of $V_\alpha$,
and that $p$ is minimal among all Klein bottles in $M(\alpha)$.
Let $P=\widehat{P}\cap M$.
We orient the boundary components of $P$ coherently on $T_0$.
Then Lemma \ref{lem:jumping} holds, because it is irrelevant to the orientability of the surfaces $P$ and $T$.

As in \cite[Section 8]{GT}, we can define two graphs $G_P$ on $\widehat{P}$ and $G_T^P$ on $\widehat{T}$
from the arcs in $P\cap T$.
We abbreviate $G_T^P$ to $G_T$.
Although $\widehat{P}$ is non-orientable,
we can assign an orientation to each vertex of $G_P$ from the orientation of $\partial u_i$.
Let $e$ be an edge of $G_P$.
If $e$ is a loop based at $u$, then
$e$ is \textit{positive\/} if a regular neighborhood $N(u\,\cup\,e)$ on $\widehat{P}$ is an annulus, \textit{negative\/} otherwise.
Assume that $e$ connects distinct vertices $u_i$ and $u_j$.
Then $N(u_i\,\cup\, e\,\cup\,u_j)$ is a disk.
Then $e$ is \textit{positive\/} if we can give an orientation to the disk $N(u_i\,\cup\, e\,\cup\, u_j)$ so that
the induced orientations on $u_i$ and $u_j$ are compatible with the original orientations of $u_i$ and $u_j$ simultaneously.
Otherwise, $e$ is \textit{negative\/}.
Then the parity rule (Lemma \ref{lem:common}(1)) still holds without change.
Also, Lemma \ref{lem:common}(2) is true.

\begin{lemma}\label{lem:GP}
$G_P$ satisfies the following.
\begin{itemize}
\item[(1)] If $t\ge 4$, then any family of parallel edges contains at most $2t$ edges.
\item[(2)] If $t\ge 3$, then any family of mutually parallel positive edges contains at most $t/2+2$ edges.
Moreover, if it contains $t/2+2$ edges, then $t\equiv 0\pmod{4}$, and $M(\beta)$ contains a Klein bottle.
\item[(3)] If $t\ge 2$ and there is a positive edge, then any family of mutually parallel negative edges contains
at most $t$ edges.
\end{itemize}
\end{lemma}

\begin{proof}
(1) is the same as Lemma \ref{lem:GS}(4).
For (2), see \cite[Lemmas 2.4(1) and 8.2(1)]{GT}.
(3) is \cite[Lemma 8.2(2)]{GT}.
\end{proof}

\begin{figure}[tb]
\includegraphics*[scale=0.45]{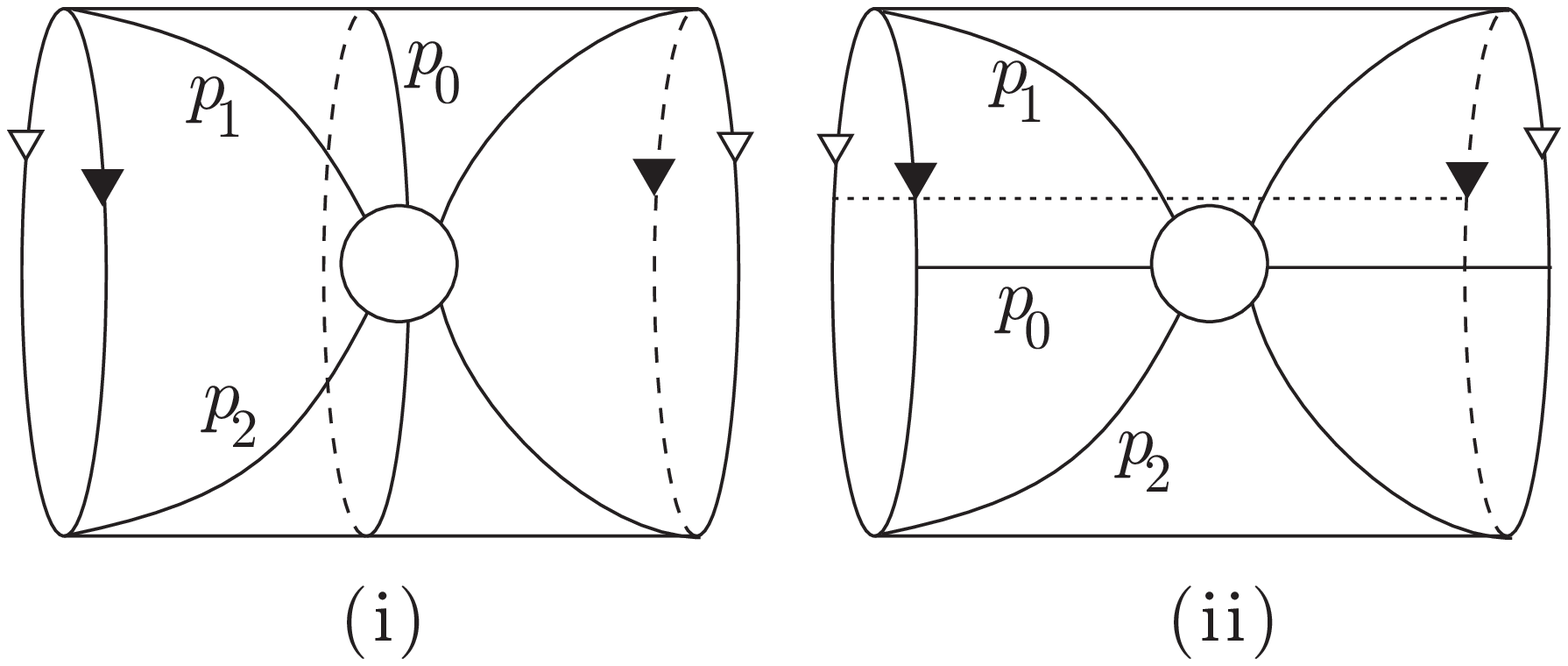}
\caption{}\label{fig:p1}
\end{figure}

When $p=1$, as in \cite[Section 8]{GT},
$G_P\cong H(p_0,p_1,p_2)$ or $H'(p_0,p_1,p_2)$, which are shown in Figure \ref{fig:p1}.
Each $p_i$ denotes the number of edges in the family of mutually parallel edges.

\section{The case where $t\le 2$}

\begin{proposition}\label{pro:t1kb}
If $t=1$, then $M$ is a $\mathbb{Q}$-homology solid torus.
\end{proposition}

\begin{proof}
Suppose $t=1$.
In the proof of \cite[Proposition 8.5]{GT}, we showed that $p=2$ and
the pair $\{G_P,G_T\}$ is uniquely determined as shown in Figure \ref{fig:kbt1}.
In fact, there is the unique edge correspondence between the edges of $G_P$ and $G_T$, and
the jumping number is one.

\begin{figure}[tb]
\includegraphics*[scale=0.6]{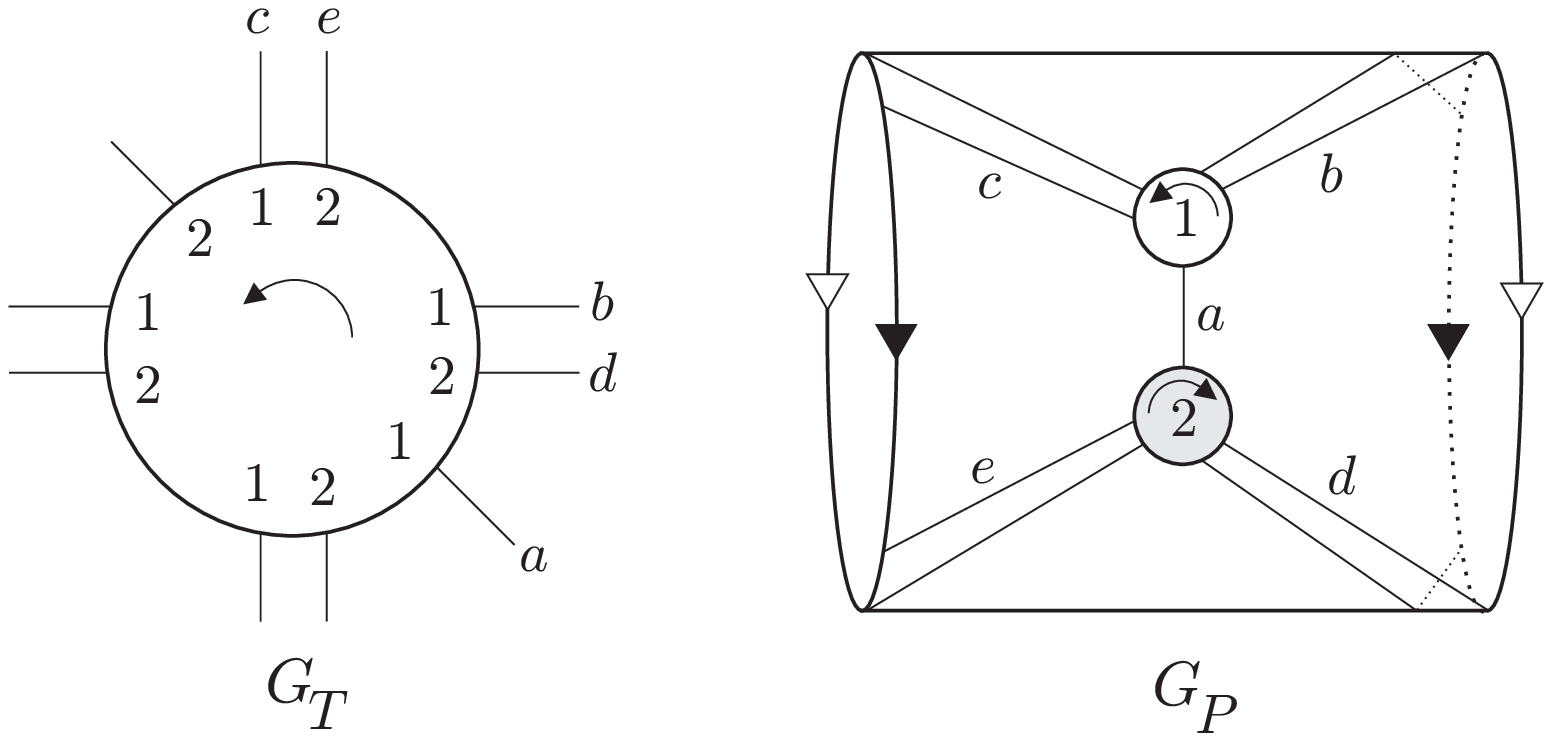}
\caption{}\label{fig:kbt1}
\end{figure}

\begin{figure}[tb]
\includegraphics*[scale=0.4]{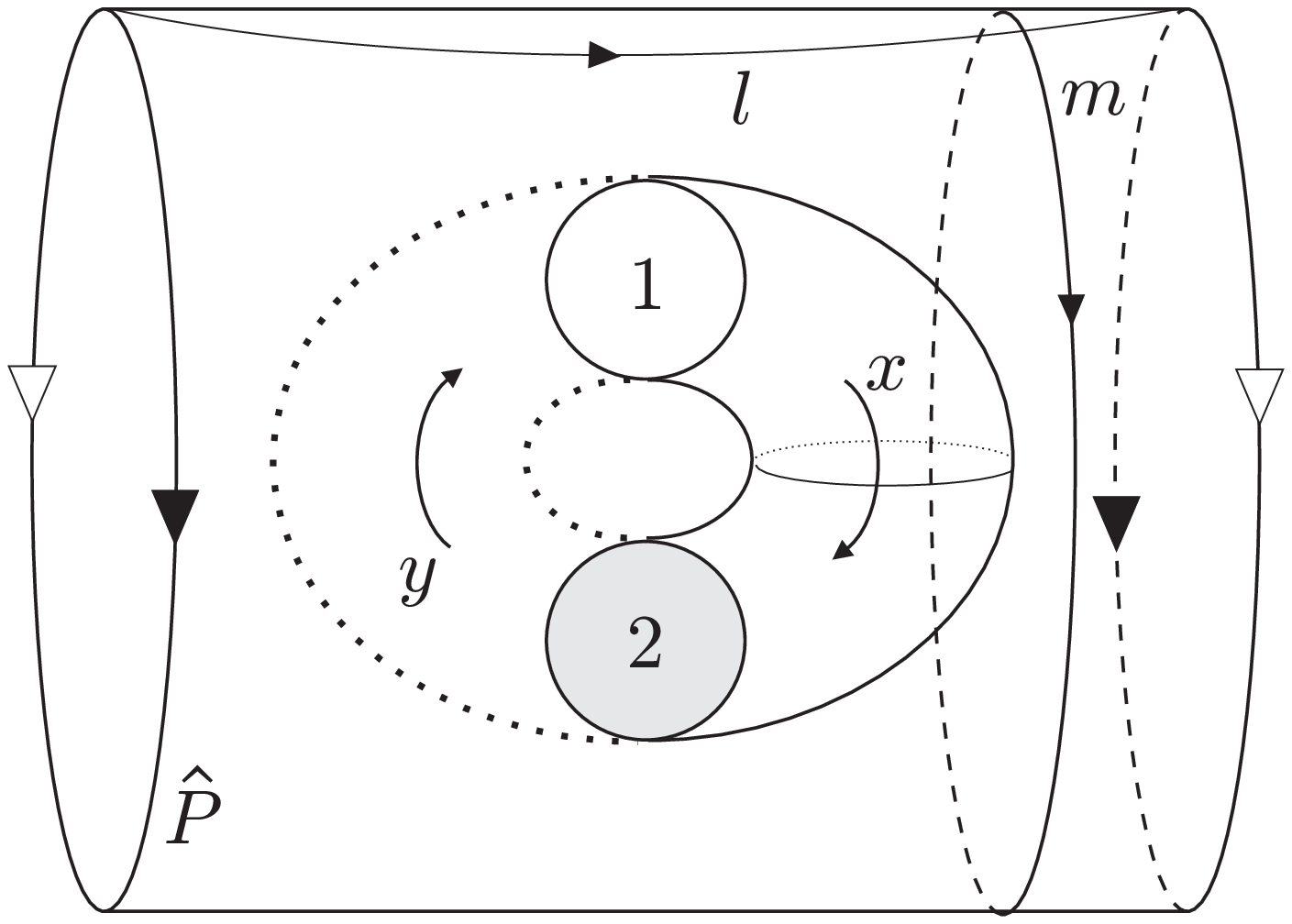}
\caption{}\label{fig:generator}
\end{figure}

We calculate $H_1(M(\alpha))$.
Let $\ell$, $m$, $x$ and $y$ be the generators of $H_1(\widehat{P}\cup V_{\alpha})$ as shown in Figure \ref{fig:generator},
where $x$ is represented by the union of the core of half of $V_\alpha$ and the edge $a$, and $y$ is similar.
Then $H_1(\widehat{P}\cup V_{\alpha})=\langle \ell,m,x ,y\ |\ 2m=0\rangle$.
Let $f$ be one of the bigons and $g_1$, $g_2$ be the $3$-gons in $G_T$.
It is easy to see that $M(\alpha)=N(\widehat{P}\cup V_\alpha\cup f_1\cup g_1\cup g_2)\cup \text{(a $3$-ball)}$. 
Hence $H_1(M(\alpha))$ has a presentation
\[
\langle \ell,m,x,y\ |\ 2m=0,x+y+2\ell+m=0,x+2y=2\ell+m,2x+y=2\ell+m  \rangle
\]
We see that $H_1(M(\alpha))=\mathbb{Z}_{20}$.
Thus $M$ is a $\mathbb{Q}$-homology solid torus.
\end{proof}

Next, we assume $t=2$.
Recall that two vertices of $G_T$ have opposite signs \cite[Lemma 9.1]{GT}.

\begin{lemma}
$p\le 2$.
\end{lemma}

\begin{proof}
This follows from \cite[Lemmas 9.6, 9.8 and 9.11]{GT}.
\end{proof}

\begin{lemma}\label{lem:t2p1}
If $p=1$, then $M$ is a $\mathbb{Q}$-homology solid torus.
\end{lemma}

\begin{proof}
By \cite[Proposition 8.7]{GT}, there are only two possibilities for $G_P$: $H(3,1,1)$ and $H'(3,2,0)$.
If $G_P=H(3,1,1)$, then $G_T=G(1,1,1,1,0)$.
The edge correspondence is shown in Figure \ref{fig:t2p1-0}.
We see that the jumping number is one.
Let $f_1$ and $f_2$ be the $3$-gons in $G_T$. 
Then we see that $M(\alpha)=N(\widehat{P}\cup V_\alpha\cup f_1\cup f_2)\cup \text{(a $3$-ball)}$, and
$H_1(M(\alpha))=\mathbb{Z}_{4}$.  
Hence $M$ is a $\mathbb{Q}$-homology solid torus.
We omit the detail.

\begin{figure}[tb]
\includegraphics*[scale=0.6]{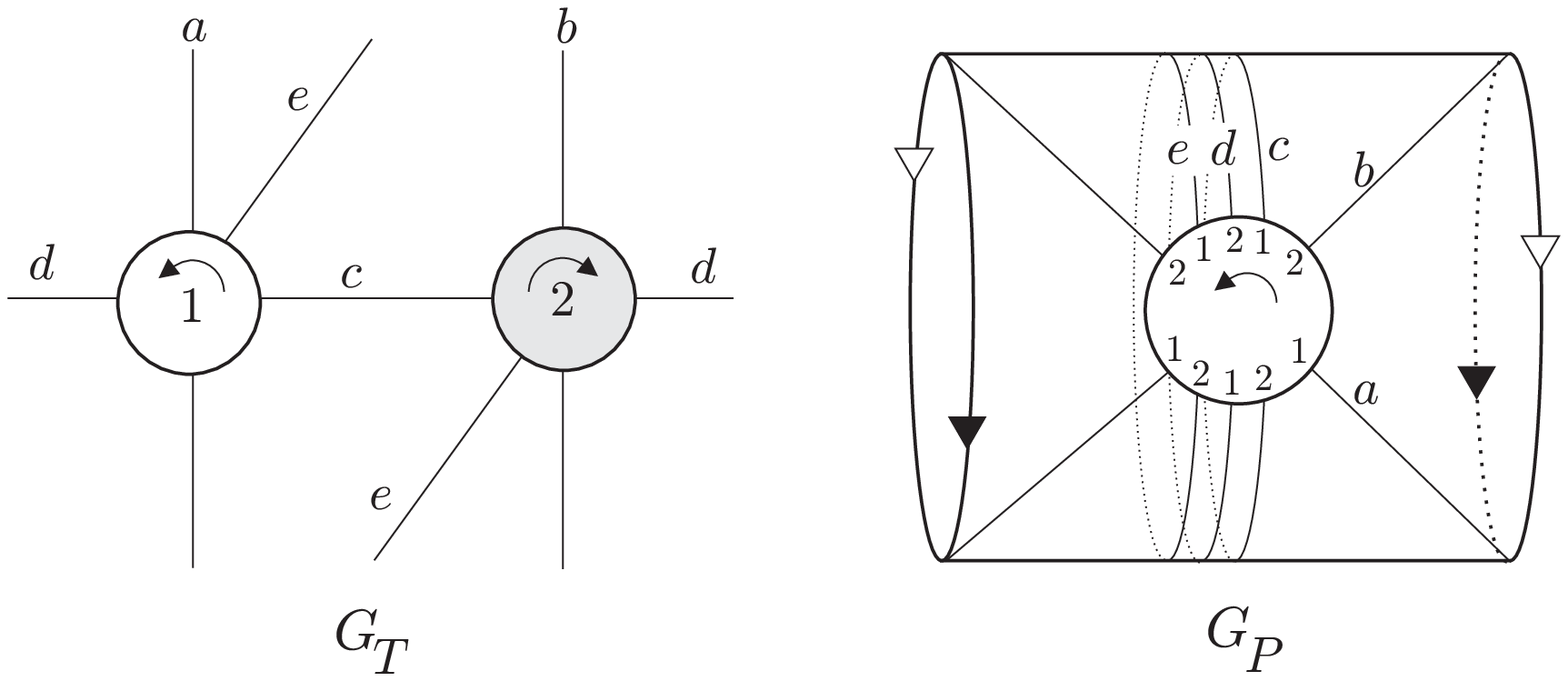}
\caption{}\label{fig:t2p1-0}
\end{figure}

If $G_P=H'(3,2,0)$, then $G_T=G(1,1,1,1,0)$ again.
By looking at the endpoints of a loop at $v_1$, the jumping number is two.
The correspondence between the edges of $G_P$ and $G_T$ is shown in Figure \ref{fig:t2p1}.

\begin{figure}[tb]
\includegraphics*[scale=0.6]{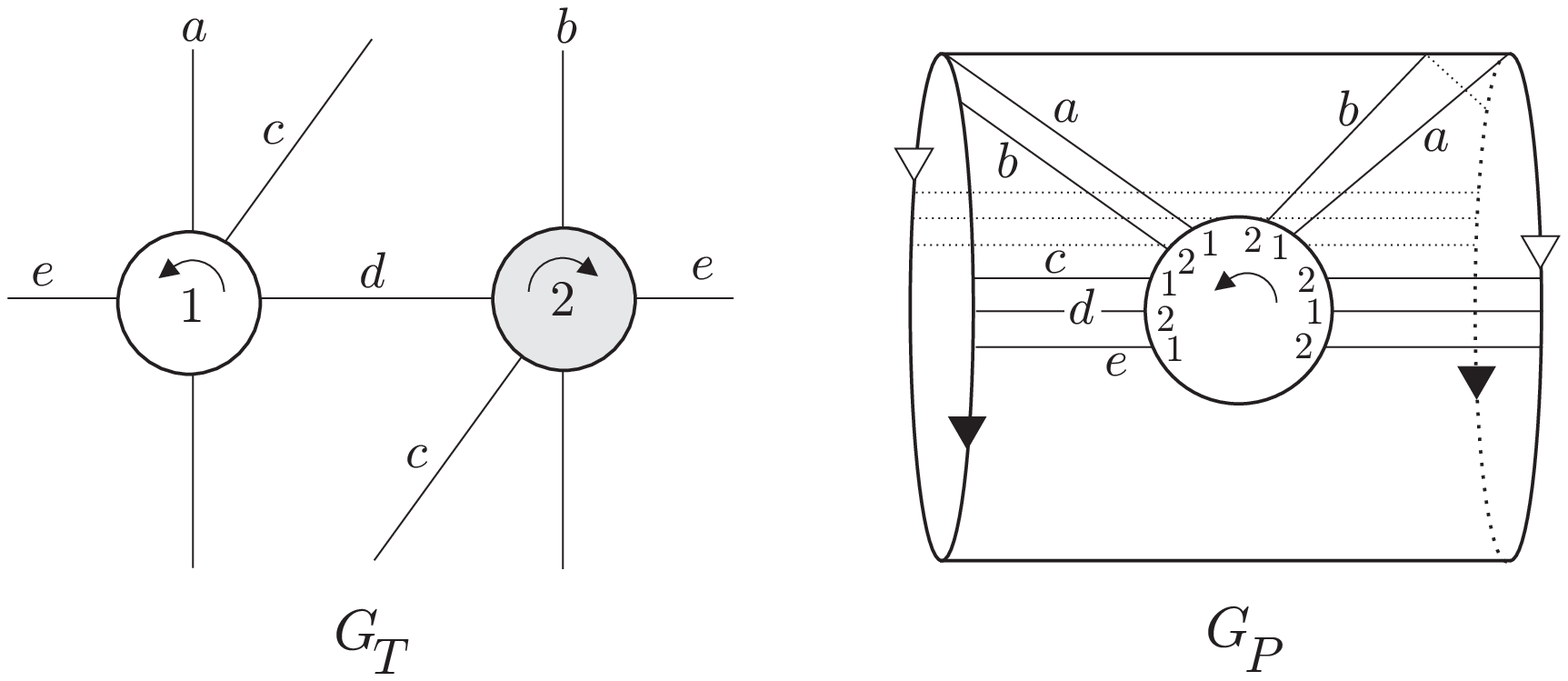}
\caption{}\label{fig:t2p1}
\end{figure}

We will calculate $H_1(M(\beta))$.
Let $f_1$ be the bigon bounded by $\{c,d\}$ in $G_P$.
We call the side of $\widehat{T}$ containing $f_1$, $\mathcal{B}$, the other side $\mathcal{W}$.
Let $g_1$ and $g_2$ be the bigon bounded by $\{d,e\}$ and the $3$-gon bounded by $\{a,b,c\}$ in $G_P$.
If we use the generators $\ell,m,x,y$ of $H_1(\widehat{T}\cup V_\beta)$ as in Figure \ref{fig:t1s2-hom},
$[\partial f_1]=2x+m$, $[\partial g_1]=2y-\ell$, and $[\partial g_2]=-y+3m$.
Let $X=N(\widehat{T}\cup V_\beta\cup f_1\cup g_1\cup g_2)$.
Then $\partial X$ has a $2$-sphere component in $\mathcal{W}$ and a torus component in $\mathcal{B}$.
Since $M(\beta)$ is irreducible, the $2$-sphere component bounds a $3$-ball in $\mathcal{W}$.
The torus component also bounds a solid torus $J$, because $M$ is hyperbolic.
Notice that $H_1(\partial J)=\langle m,\ell+x\rangle$.
Thus if the meridian of $J$ is $r(\ell+x)+sm$ for some $r,s$, then
$H_1(M(\beta))$ has a presentation
\[
\langle \ell,m,x,y\ | \ 2x+m=0,2y-\ell=0,-y+3m=0,r(\ell+x)+sm=0 \rangle,
\]
which gives $\langle x\ |\ (11r+2s)x=0\rangle$.
Unless $11r+2s=0$, $H_1(M(\beta))$ is a torsion, and $M$ is a $\mathbb{Q}$-homology solid torus.
If $11r+2s=0$, then $(r,s)=(\pm 2,\mp 11)$.
As it is well known (\cite[Theorem 4.1]{GL5}), $\mathcal{B}$ admits a Seifert fibration over the disk with two exceptional fibers, one of which
has index two.
Moreover, a regular fiber represents $m$ on $\widehat{T}$.
Hence if $(r,s)=(\pm 2,\mp 11)$, then the regular fiber intersects the meridian of $J$ just twice.
This implies that another exceptional fiber of $\mathcal{B}$, which is a core of $J$, has index two, and hence
$\mathcal{B}$ contains a Klein bottle as in the proof of Lemma \ref{lem:t2opposite}.
Thus $M$ is a $\mathbb{Q}$-homology solid torus by Lemma \ref{lem:beta-noklein}.
\end{proof}


\begin{lemma}\label{lem:t2p2}
If $p=2$, then $M$ is a $\mathbb{Q}$-homology solid torus.
\end{lemma}

\begin{proof}
By \cite[9.1 and 9.2]{GT}, there are only two possibilities for $G_T$: $G(1,2,2,2,2)$ and $G(2,2,1,2,1)$.

Suppose $G_T=G(1,2,2,2,2)$.
Then the associated permutation to each pair of parallel negative edges is $(12)$ by \cite[Lemma 9.3]{GT}.
Two edges in each pair form an essential orientation-preserving cycle on $\widehat{P}$ \cite[Lemma 2.3]{Go2}.
Hence $G_P$ has four mutually parallel positive edges.
Among them, there are two bigons lying in the same side of $\widehat{T}$.
If $M(\beta)$ does not contain a Klein bottle, then these bigons have the same edge class labels by Lemma \ref{lem:atmost3}(1).
This contradicts Lemma \ref{lem:common}(2).
If $M(\beta)$ contains a Klein bottle, then $M$ is a $\mathbb{Q}$-homology solid torus by Lemma \ref{lem:beta-noklein}.

Suppose $G_T=G(2,2,1,2,1)$.
Then the associated permutation to two pairs of parallel negative edges is $(12)$ by \cite[Lemma 9.9]{GT}.
Hence each vertex of $G_P$ is incident to one positive loop and two negative loops, and there are four positive edges
between $u_1$ and $u_2$.
We see that the jumping number is two, and $G_P$ is uniquely determined as shown in Figure \ref{fig:22121}.

\begin{figure}[tb]
\includegraphics*[scale=0.65]{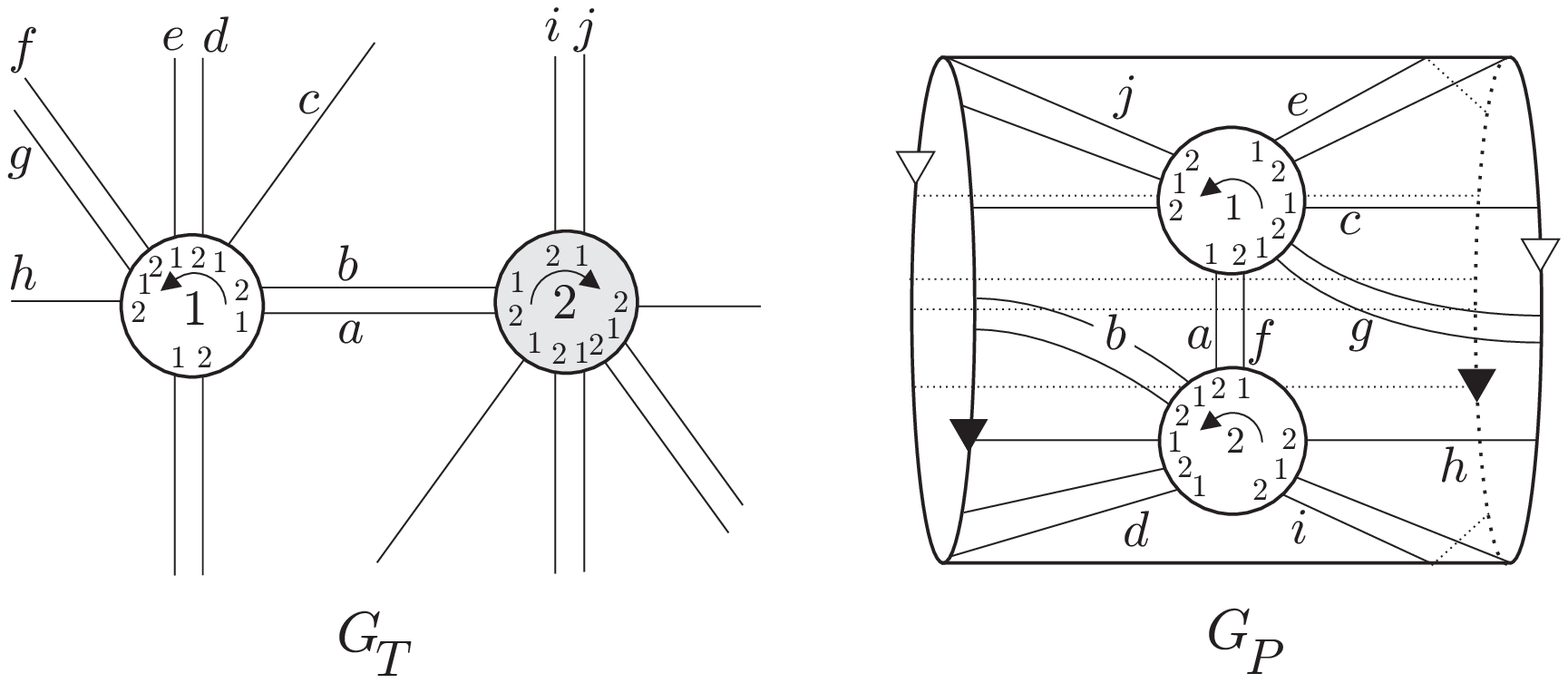}
\caption{}\label{fig:22121}
\end{figure}

Let $f_1$ and $f_2$ be the bigons bounded by $\{a,b\}$ and $\{d,e\}$, respectively in $G_T$.
Let $g$ be the $3$-gon bounded by $\{a,c,i\}$.
By using the generators of $H_1(\widehat{P}\cup V_\alpha)$ as in Figure \ref{fig:generator},
$[\partial f_1]=x-y-2\ell$, $[\partial f_2]=x+y+m$, and $[\partial g]=-x+5\ell+m$.
Thus $H_1(M(\alpha))$ has a presentation
\[
\langle \ell,m,x,y\ |\ 2m=0,x-y-2\ell=0,x+y+m=0,-x+5\ell+m=0\rangle,
\]
giving $H_1(M(\alpha))=\mathbb{Z}_{16}$.
Thus $M$ is a $\mathbb{Q}$-homology solid torus.
\end{proof}

\section{The case where $t\ge 3$}

Finally, we assume $t\ge 3$.

\begin{lemma}\label{lem:p1t3}
$p=1$.
\end{lemma}

\begin{proof}
The possibility $p\ge 3$ was ruled out in Subsection 11.1 of \cite{GT}.
Also, Subsection 11.2 of \cite{GT} eliminates the possibility $p=2$, but
we will give a detail of the remark after \cite[Lemma 11.20]{GT}.

Now, $G_P$ has no positive edges.
Each vertex is incident to a family of $n$ mutually parallel negative loops, where $n\ge 2t$.

First, assume $t\ge 4$.
Then $n=2t$ by Lemma \ref{lem:GP}(1).
Hence $G_P$ contains just $t$ non-loop negative edges.
Let $\sigma$ be the associated permutation to the family of negative loops at $u_1$.
Then $\sigma$ has a single orbit as in Lemma \ref{lem:GS}(3).
If the edge endpoints of non-loop edges are successive around $u_1$, then
$\sigma$ would be the identity.
Hence the edge endpoints of non-loop edges are not successive around $u_1$ (and so $u_2$).
Then the reduced graph of $G_P$ has the form as shown in Figure \ref{fig:p2kb}.

\begin{figure}[tb]
\includegraphics*[scale=0.6]{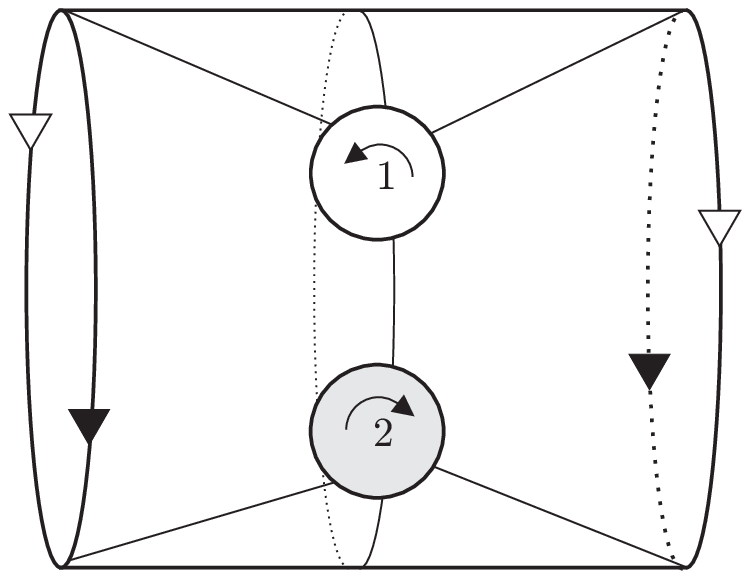}
\caption{}\label{fig:p2kb}
\end{figure}

Let $A_1,\dots,A_t,B_1,\dots,B_t,C_1,\dots,C_t$ be the edges of $G_P$ incident to $u_1$ as in Figure \ref{fig:p2kb-label}.
Here, $\sigma(i)\equiv i+a \pmod{t}$.
By symmetry, we may assume that $1\le a\le t/2$.
If $a=t/2$, then $\sigma^2$ would be the identity.
This contradicts that $\sigma$ has a single orbit, since $t\ge 4$.
Hence $2a+1\le t$.

\begin{figure}[tb]
\includegraphics*[scale=0.45]{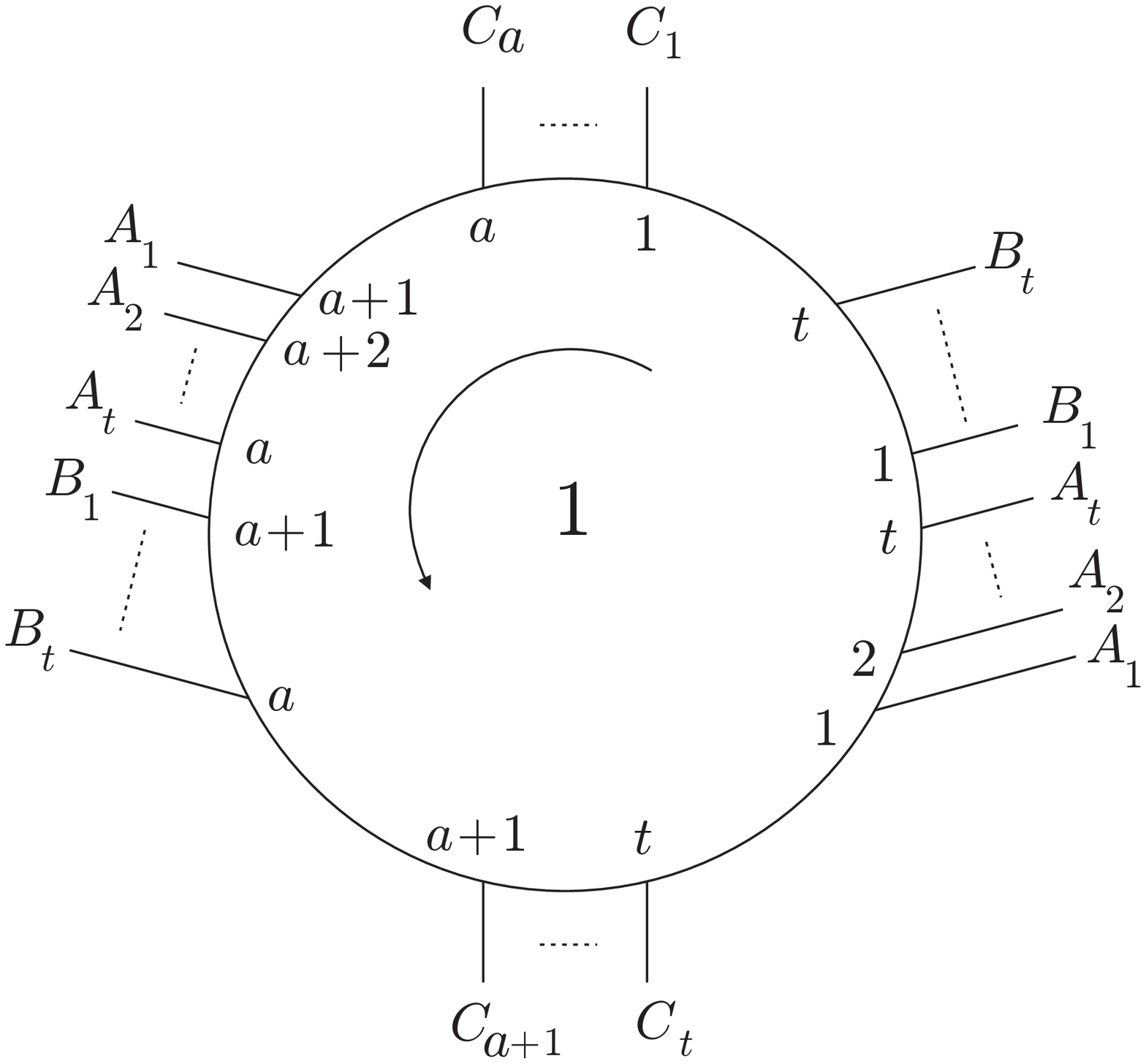}
\caption{}\label{fig:p2kb-label}
\end{figure}

Let $S$ be the cycle formed by the edges $A_1,\dots,A_t$ on $\widehat{T}$.
It is an essential orientation-preserving loop there.

\begin{claim}\label{cl:jump1}
The jumping number is not one.
\end{claim}

\begin{proof}[Proof of Claim \ref{cl:jump1}]
Assume that the jumping number is one.
Then $S\cup B_1\cup B_{a+1}$ is as shown in Figure \ref{fig:twoiso}(i).

\begin{figure}[tb]
\includegraphics*[scale=0.4]{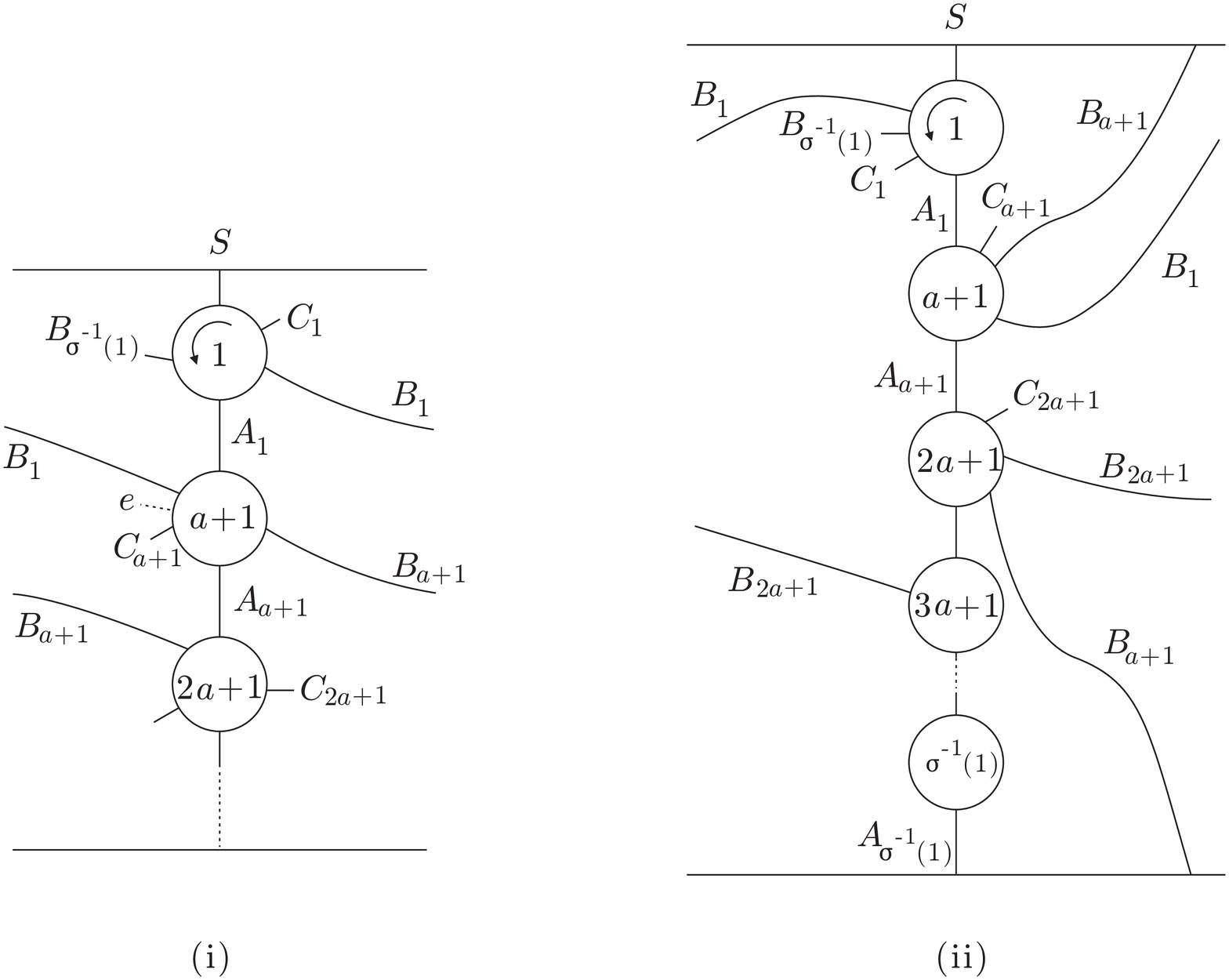}
\caption{}\label{fig:twoiso}
\end{figure}

If $C_{a+1}$ goes to $v_{2a+1}$, then it is parallel to $A_{a+1}$.
Then there is no edge incident to the edge endpoint with label $2$ at $v_{a+1}$ between $A_{a+1}$ and $C_{a+1}$.
By a similar argument, $C_{a+1}$ cannot go to $v_{1}$.
Thus $C_{a+1}$ goes to $v_{a+1}$.
Then we see that all $C_i$'s have the same label at its endpoints.
Moreover, the family of mutually parallel negative edges at $u_2$ has the same permutation $\sigma$.

Let $e$ be the edge which is incident to $v_{a+1}$ between $B_1$ and $C_{a+1}$ (see Figure \ref{fig:twoiso}(i)).
Then it goes to $v_1$.
Let $p$ and $q$ be the edge endpoints of $C_{a+1}$ and $e$ at $v_{a+1}$ with label $2$, respectively.
Then $p$ and $q$ are not successive among five points of $\partial v_{a+1}\cap \partial u_2$.
In fact, they appear in the order $p,*,q,*,*$ on $\partial v_{a+1}$ along its orientation.
Hence we see that $e$ goes to $v_{2a+1}$ by examining the five occurrences of labels $a+1$ around $u_2$, a contradiction.
\end{proof}

\begin{claim}\label{cl:jump2}
The jumping number is not two.
\end{claim}

\begin{proof}[Proof of Claim \ref{cl:jump2}]
Assume that the jumping number is two.
The arrangements of $S$, $B_1$, $B_{a+1}$ are shown in Figure \ref{fig:twoiso}(ii).
If $t=4$, then $a=1$, and we cannot locate $B_3$.
Otherwise, $B_{2a+1}$ should be as in Figure \ref{fig:twoiso}(ii).
Then we cannot locate $B_{\sigma^{-1}(1)}$ there.
\end{proof}

Thus we have eliminated the case $t\ge 4$.

Next, assume $t=3$.
Let $k$ be the number of non-loop edges in $G_P$.
Since $k+2n=5t=15$ and $n\ge 2t=6$, $k\le 3$.
In fact, $k=1$ or $3$, because $k$ is odd.

\begin{claim}\label{cl:k1}
$k\ne 1$.
\end{claim}

\begin{proof}[Proof of Claim \ref{cl:k1}]
Let $k=1$.  Then $n=7$.
Let $A_1,A_2,A_3$, $B_1,B_2,B_3$, $C_1$ be the negative loops at $u_1$ numbered successively.
We may assume that $A_1, B_1, C_1$ have label $1$ at one end of the family.
Thus the associated permutation to the family is $(123)$.
The edges $A_1,A_2,A_3$ form an essential cycle $S$ on $\widehat{T}$.
The arrangements of $B_1$ and $C_1$ are shown in Figure \ref{fig:k1}(i) and (ii), according to the jumping number.
In any case, $B_2$ cannot be placed.
\end{proof}

\begin{figure}[tb]
\includegraphics*[scale=0.4]{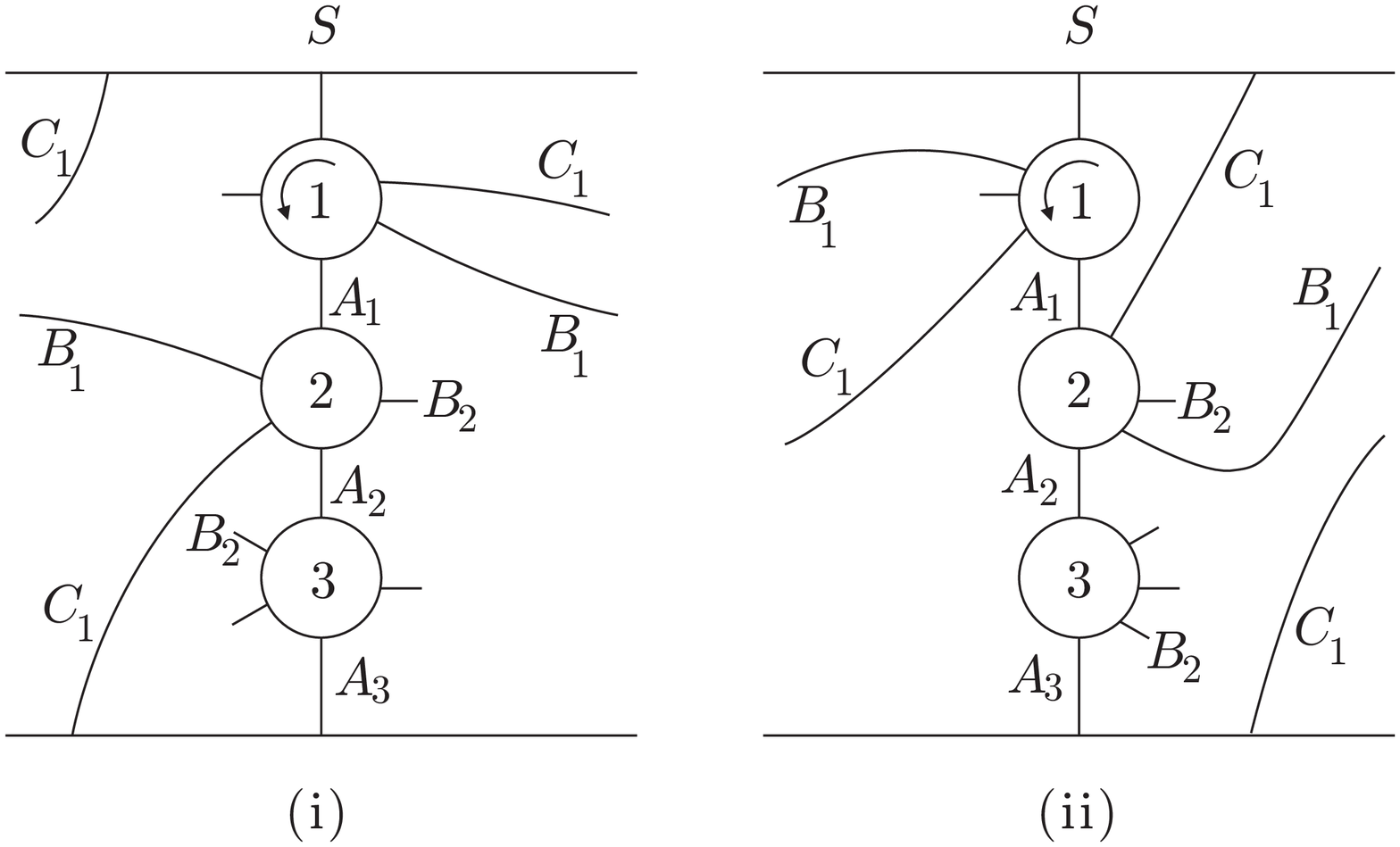}
\caption{}\label{fig:k1}
\end{figure}

\begin{claim}\label{cl:k3}
$k\ne 3$.
\end{claim}

\begin{proof}[Proof of Claim \ref{cl:k3}]
Let $k=3$.  Then $n=2t=6$.
Hence the edge endpoints of non-loop negative edges are not successive as before.
The proof of Claim \ref{cl:jump1} works here without any change, hence the jumping number is not one.
Let $A_1,A_2,A_3,B_1,B_2,B_3$ be the negative loops at $u_1$ numbered successively 
and let $C_1,C_2,C_3$ be the non-loop negative edges as in Figure \ref{fig:k3}.

\begin{figure}[tb]
\includegraphics*[scale=0.4]{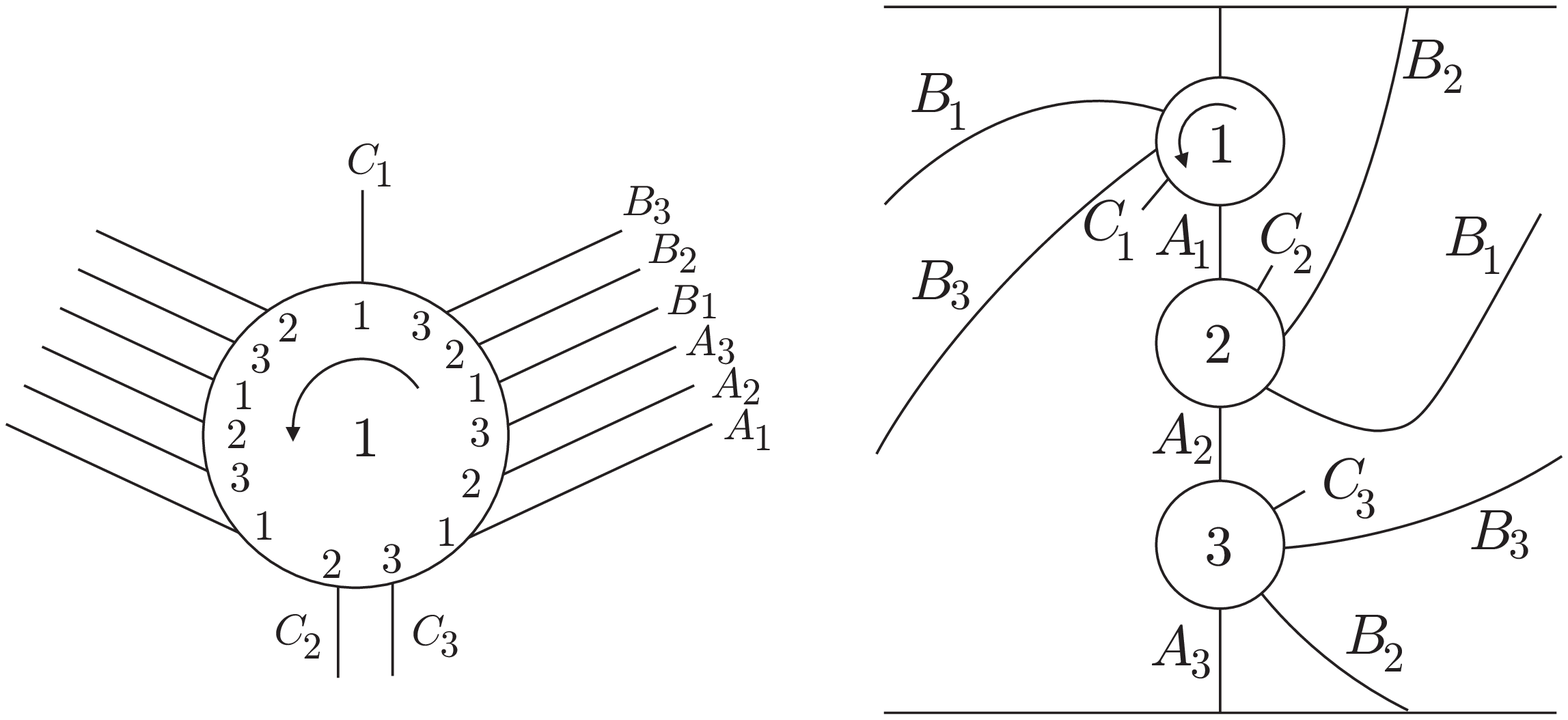}
\caption{}\label{fig:k3}
\end{figure}

Consider $C_2$.
If it goes to $v_3$, then it is parallel to $B_2$.
However, $B_2$ has label $1$ at both endpoints, $C_2$ has label $1$ at $v_2$ but label $2$ at $v_3$.
Clearly, this is impossible.
Hence $C_2$ goes to $v_1$, so it is parallel to $A_1$.
This is also impossible by the same reason.
\end{proof}

This completes the proof of Lemma \ref{lem:p1t3}
\end{proof}


Since the vertex of $G_P$ has degree $5t$, $t$ is even.
Thus $t\ge 4$.
We distinguish two cases.

\subsection{$G_P$ has no positive loops}

In this case, $G_P=H(0,p_1,p_2)$ where $p_1+p_2=5t/2$.
Let $U$ and $V$ be the families of mutually parallel negative loops with $p_1$ and $p_2$ edges, respectively.
Without loss of generality, we may assume $p_1>t$.
Then the associated permutation $\sigma$ to $U$ has a single orbit.
Thus all vertices of $G_T$ have the same sign.
Moreover, we see that $p_1=t+h$ for some odd number $h$ with $\gcd (t,h)=1$.
Hence, $p_1$ or $p_2>t+1$, unless $t=4$ and $p_1=p_2=t+1$.

\begin{lemma}
The case where $t=4$ and $p_1=p_2=t+1$ is impossible.
\end{lemma}

\begin{proof}
We may assume that the labels of $G_P$ are as shown in Figure \ref{fig:p1t3}.

\begin{figure}[tb]
\includegraphics*[scale=0.35]{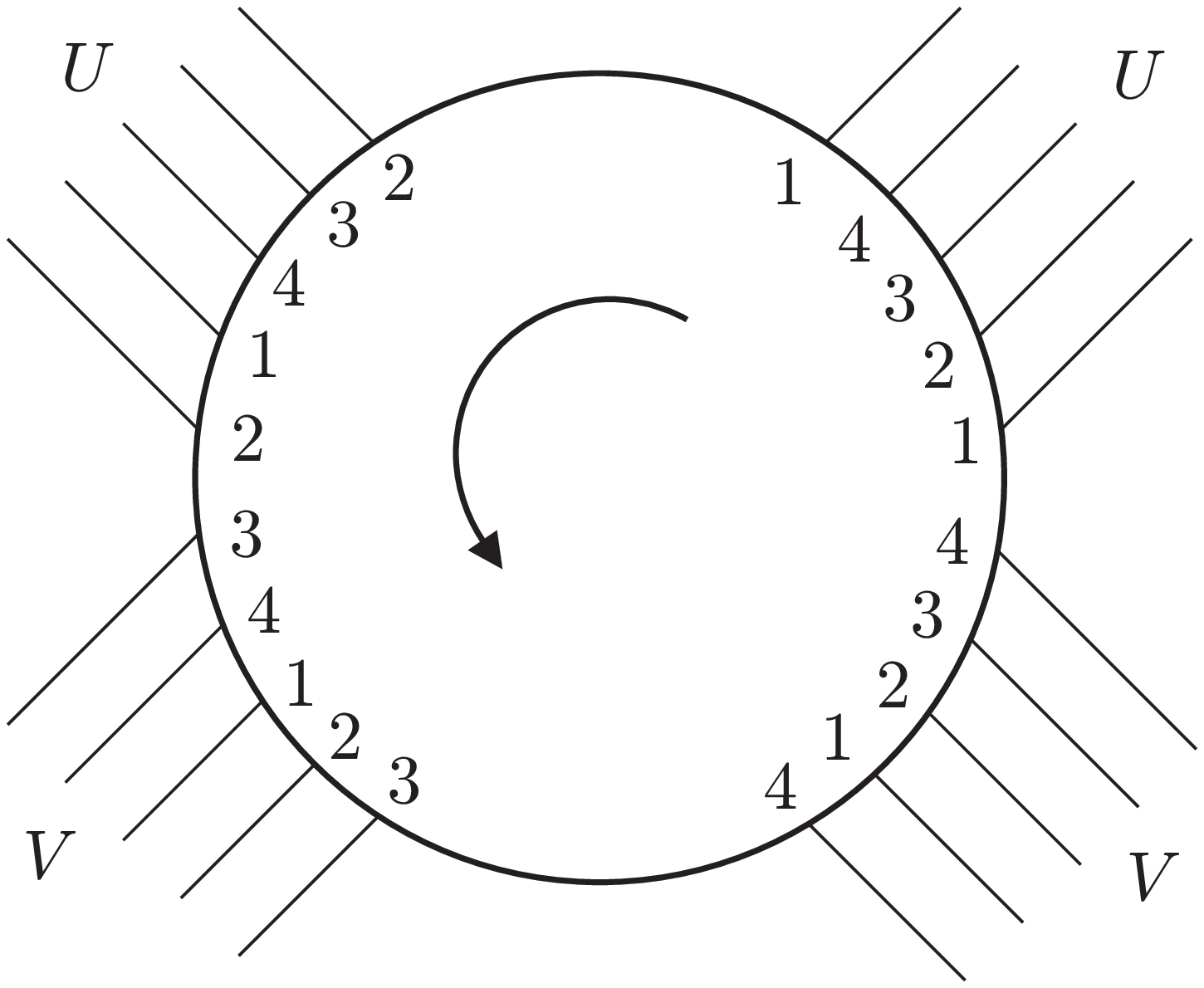}
\caption{}\label{fig:p1t3}
\end{figure}

Then $\nu(v_1,v_2)\ge 2$ and $\nu(v_3,v_4)\ge 2$, since $U$ contains two $\{1,2\}$-edges, and $V$ contains two $\{3,4\}$-edges.
(Recall that $\nu(v_i,v_j)$ denotes the number of mutually non-parallel edges between $v_i$ and $v_j$ in $G_T$.)
Thus $\nu(v_1,v_2)=2$ by \cite[Lemma 5.4(ii)]{Go2}.
Let $e$ be the $\{1,2\}$-edge in $V$.
Then $e$ is parallel to either of $\{1,2\}$-edges in $U$.
In any case, this contradicts \cite[Lemma 2.4]{Go2}.
\end{proof}

\begin{figure}[tb]
\includegraphics*[scale=0.35]{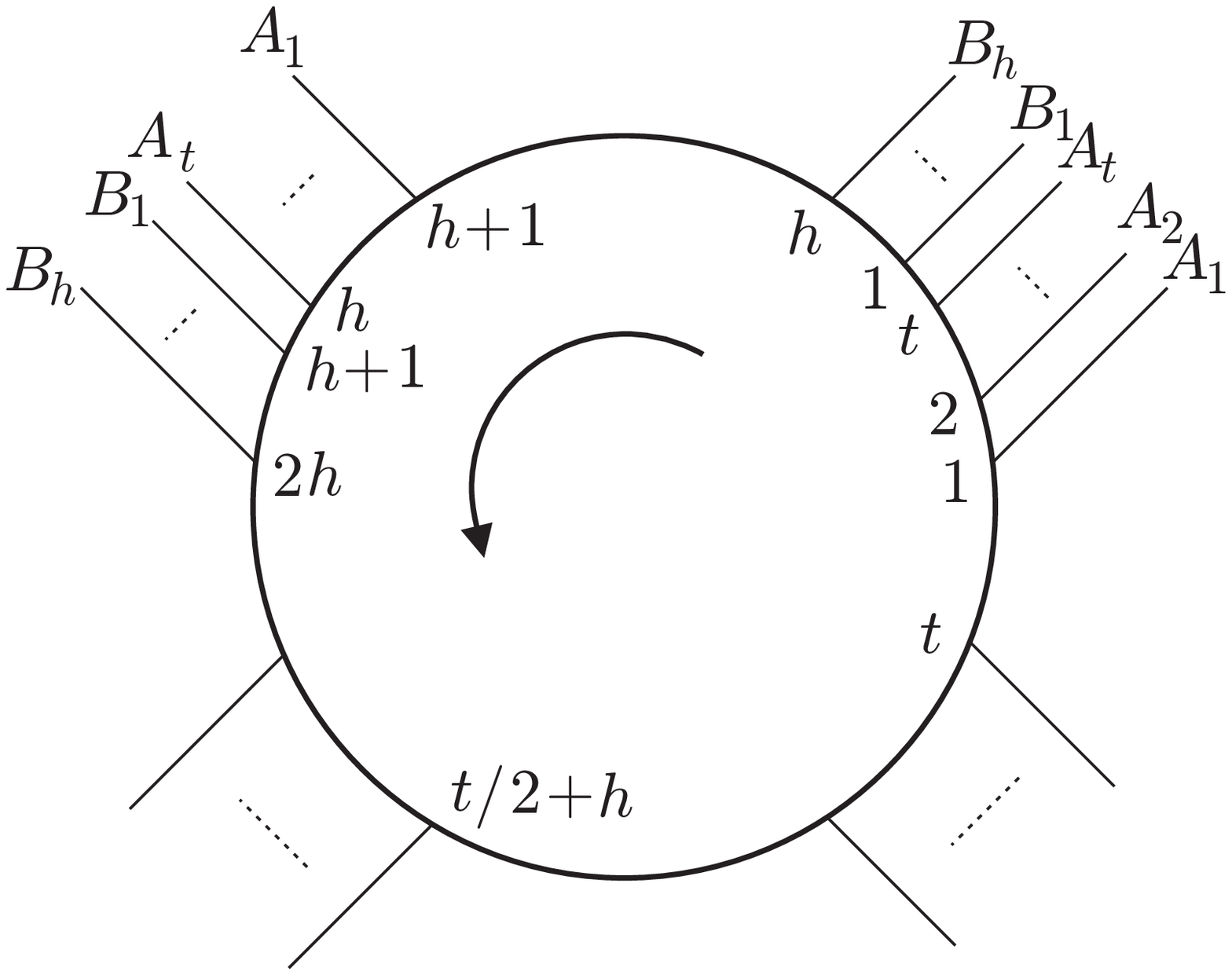}
\caption{}\label{fig:p1t3-label}
\end{figure}

Thus $p_1=t+h$ with $3\le h\le t-1$ by Lemma \ref{lem:GP}(1) and $\gcd(t,h)=1$.
We may assume that $G_P$ has the labels as shown in Figure \ref{fig:p1t3-label}.
Let $A_1,A_2,\dots,A_t,B_1,\dots,B_h$ be the edges of $U$ numbered successively such that $A_i$ and $B_i$ have label $i$ and $\sigma(i)$,
where $\sigma(i)\equiv i+h \pmod{t}$.
Since $\sigma$ has a single orbit, the edges $A_1,A_2,\dots,A_t$ form an essential cycle $S$ on $\widehat{T}$.

\begin{lemma}
$h\le t/2-1$.
\end{lemma}

\begin{figure}[tb]
\includegraphics*[scale=0.4]{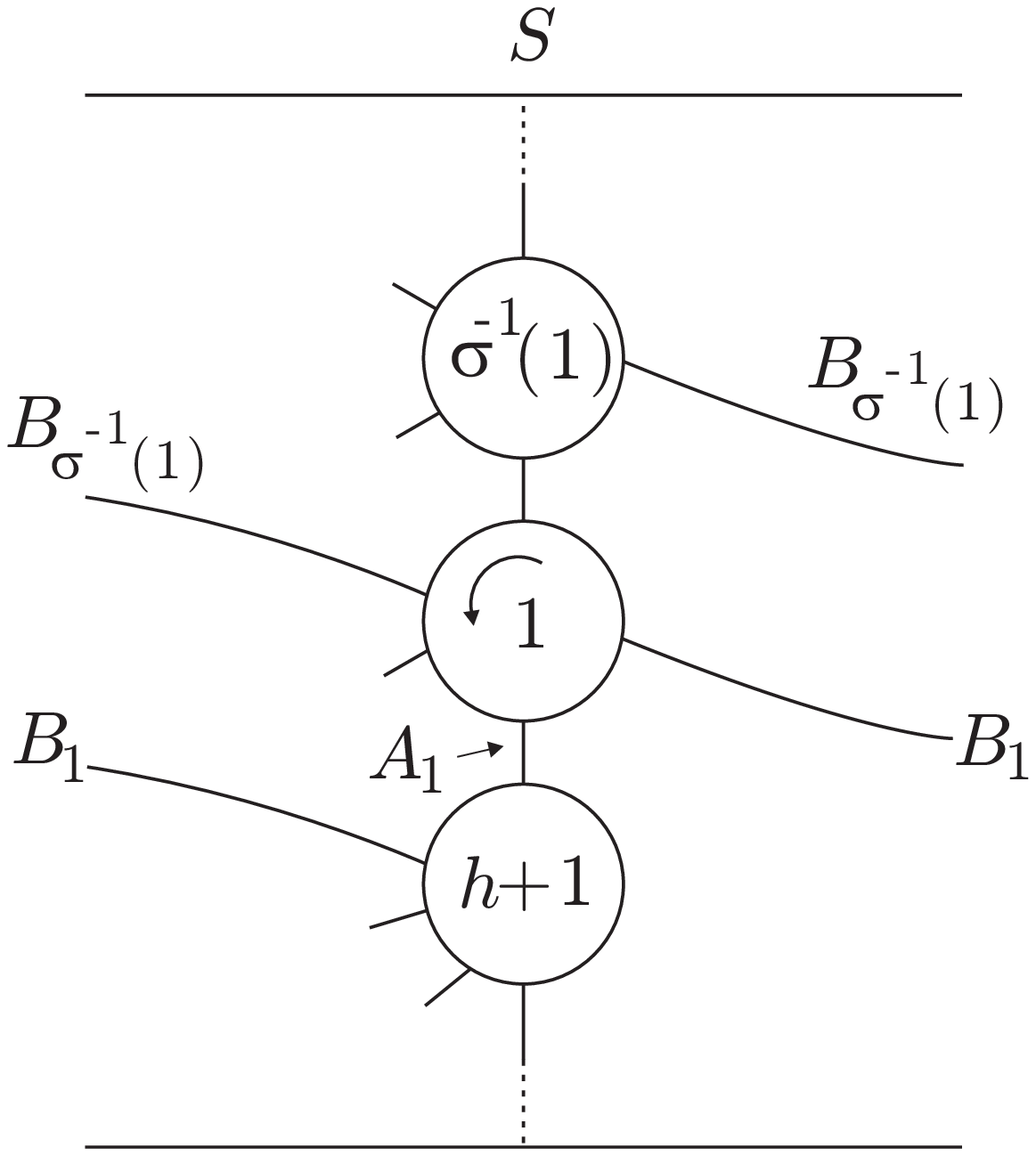}
\caption{}\label{fig:jump1}
\end{figure}

\begin{figure}[tb]
\includegraphics*[scale=0.4]{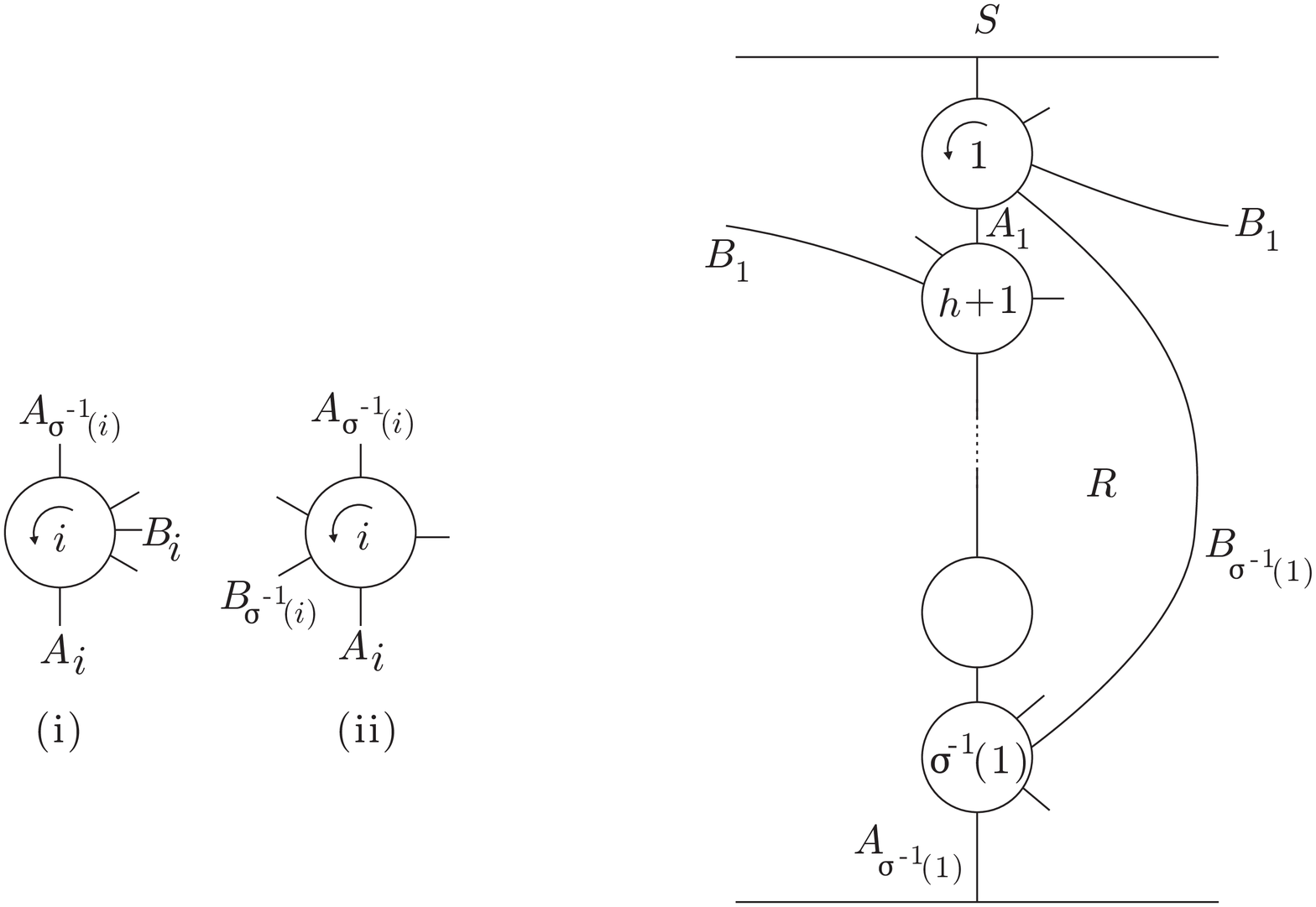}
\caption{}\label{fig:jump2-final}
\end{figure}

\begin{proof}
Suppose not.  Then $h\ge t/2+1$, since $\gcd(h,t)=1$.
Hence the label $1$ appears exactly four times at the endpoints of the edges of $U$.
Consider $S\cup B_1\cup B_{\sigma^{-1}(1)}$.

If the jumping number is one, then the situation is as in Figure \ref{fig:jump1}.
Then the edge between $A_1$ and $B_{\sigma^{-1}(1)}$ at $v_1$ cannot go to any vertex.
Thus the jumping number is two.
The situation around each vertex $v_i$ is as in Figure \ref{fig:jump2-final},
where (i) corresponds to $i=1,2,\dots,h$, and (ii) corresponds to $i=h+1,\dots,t$.
Also, Figure \ref{fig:jump2-final} shows $S\cup B_1\cup B_{\sigma^{-1}(1)}$.
Since $h\ge 3$, there is a vertex $v_j$ between $v_{h+1}$ and $v_{\sigma^{-1}(1)}$ along $S$ for $1<j\le h$.
For such $j$, $B_j$ lies in the region $R$.
Hence $B_j$ is parallel to $A_j$, a contradiction.
\end{proof}

\begin{lemma}\label{lem:positive-loop}
$G_P$ contains a positive loop.
\end{lemma}

\begin{proof}
First, suppose that the jumping number is one.
The family $V$ of $G_P$ contains a $(t/2+h,t)$-edge.
Hence no vertex on $S$ between $v_{t/2+h}$ and $v_{t}$ is incident to some $B_i$.
See Figure \ref{fig:jump1-final}.

\begin{figure}[tb]
\includegraphics*[scale=0.4]{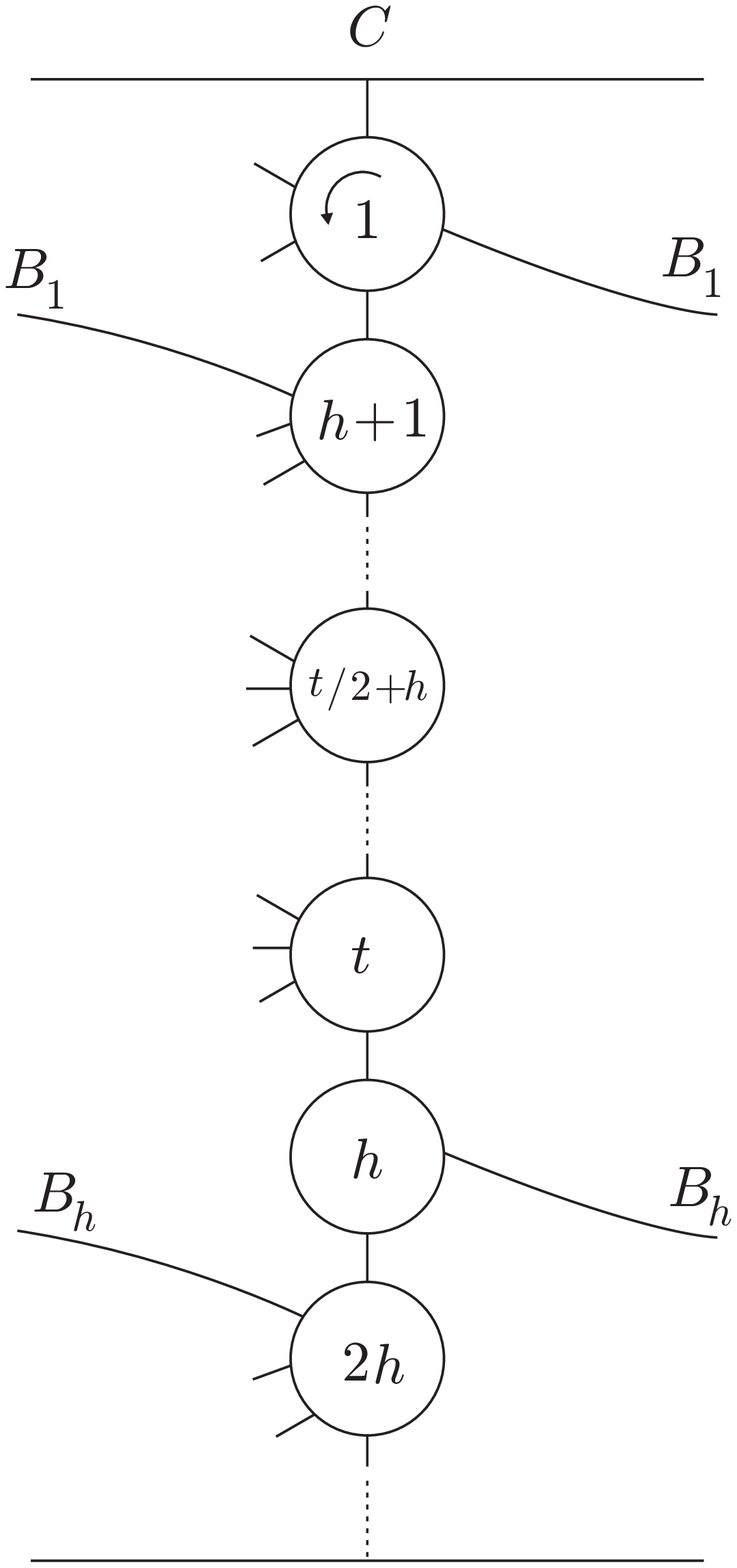}
\caption{}\label{fig:jump1-final}
\end{figure}

However, this means that $(t/2+h)+kh=t$ for some $k>0$, since $\sigma(i)\equiv i+h \pmod{t}$.
Thus $t$ is a multiple of $h$, contradicting $\gcd(h,t)=1$.

When the jumping number is two, a similar argument works.
\end{proof}

\subsection{$G_P$ has a positive loop}

Recall that $G_P\cong H(p_0,p_1,p_2)$ or $H'(p_0,p_1,p_2)$ where $p_0+p_1+p_2=5t/2$.

\begin{lemma}\label{lem:finalcases}
$p_0=t/2$, $t/2+1$ or $t/2+2$.
\end{lemma}

\begin{proof}
Since $p_0\ne 0$, $p_i\le t$ for $i=1,2$ by Lemma \ref{lem:GP}(3).
Thus $t/2\le p_0\le t/2+2$ Lemma \ref{lem:GP}(2).
\end{proof}

If $p_0=t/2+2$, then
$M(\beta)$ contains a Klein bottle by Lemma \ref{lem:GP}(2), hence $M$ is a $\mathbb{Q}$-homology solid torus by
Lemma \ref{lem:beta-noklein}.
We consider two remaining cases.

\bigskip
\textsc{Case (1)}. $p_0=t/2$.

By Lemma \ref{lem:GP}, $G_P=H(t/2,t,t)$ or $H'(t/2,t,t)$.

\begin{lemma}\label{lem:244}
If $G_P\cong H(t/2,t,t)$, then $t=4$ and $M$ is a $\mathbb{Q}$-homology solid torus.
\end{lemma}

\begin{proof}
We may assume that $t$ positive loops have labels $1,2,\dots,t/2$ at their endpoints successively.
First, we claim that $t/2$ is even.
If $t/2$ is odd, then the middle edge of the positive loops is a $((t+2)/4,(3t+2)/4)$-edge.
But a family of $t$ mutually parallel negative loops contains a negative $((t+2)/4,(3t+2)/4)$-edge, a contradiction.

Let $A$ and $B$ be the families of $t$ mutually parallel negative loops.
Then the associated permutation to $A$ (and $B$) has $t/2$ orbits of length two.
Let $H$ be the subgraph of $G_T$ spanned by $v_1$, $v_{t/2}$, $v_{t/2+1}$ and $v_t$.
If $t>4$, then $H$ has an annulus support as in Figure \ref{fig:t/2}.
A jumping number argument easily rules out this configuration.

\begin{figure}[tb]
\includegraphics*[scale=0.6]{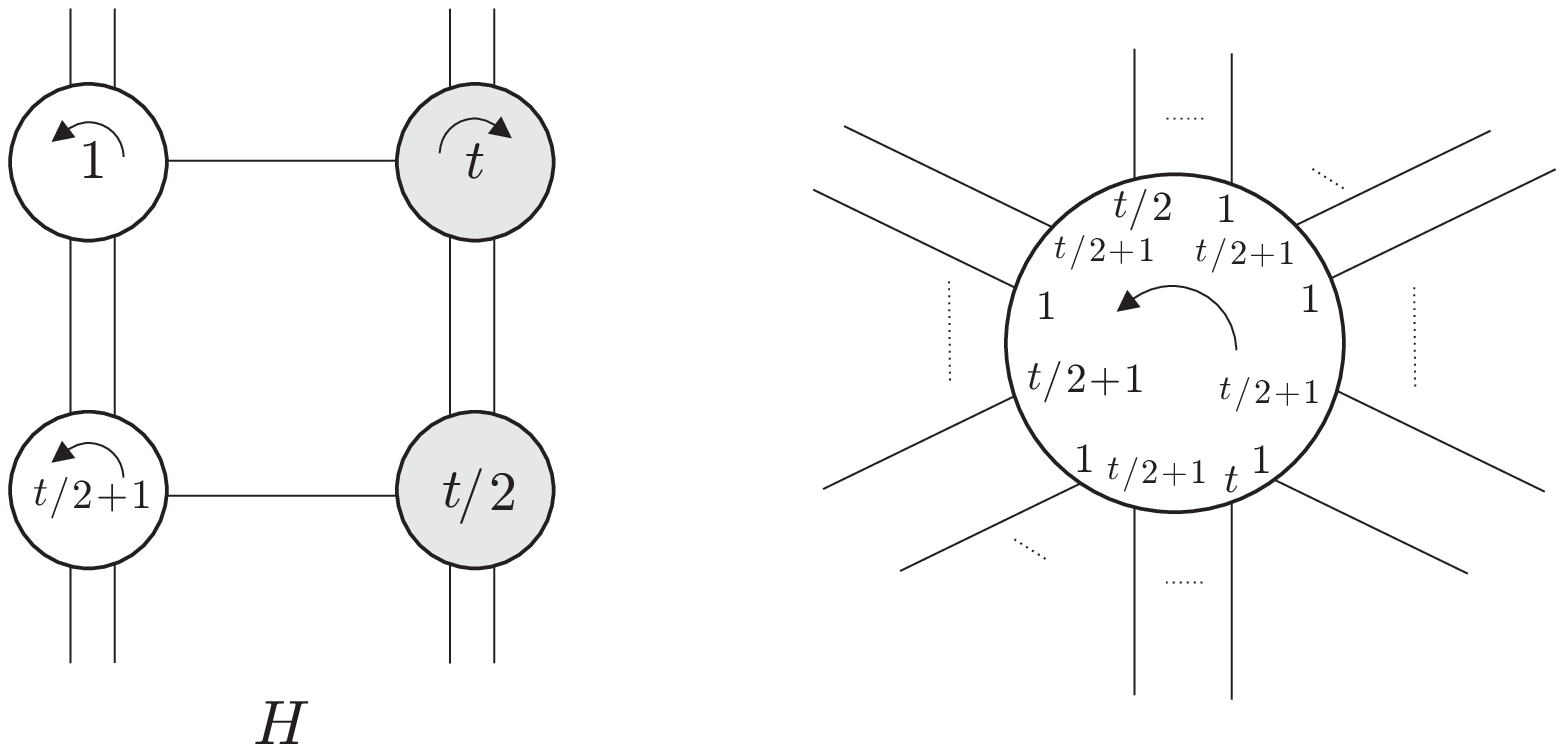}
\caption{}\label{fig:t/2}
\end{figure}

Thus $t=4$.  In fact, $G_T$ has a torus support, and we see that the jumping number is one.
Figure \ref{fig:244} shows the correspondence between the edges of $G_P$ and $G_T$.

\begin{figure}[tb]
\includegraphics*[scale=0.6]{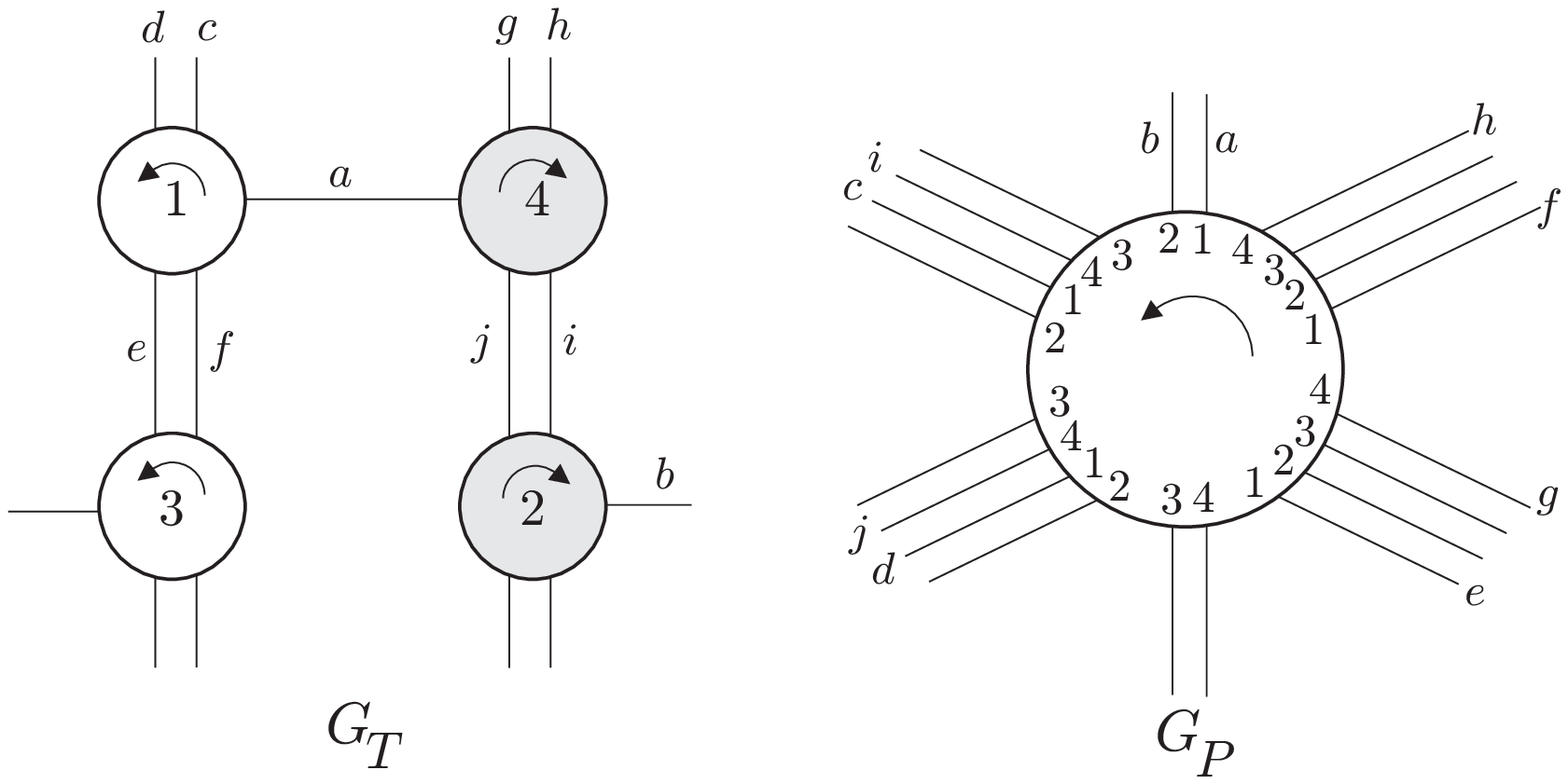}
\caption{}\label{fig:244}
\end{figure}

Let $f_1$ be the bigon bounded by $\{e,f\}$ and $f_2$ the $6$-gon bounded by $\{a,c,f,g,j\}$ in $G_T$.
Then we see that $M(\alpha)=N(\widehat{P}\cup V_\alpha\cup f_1\cup f_2)\cup \text{(a $3$-ball)}$.
It is easy to show $H_1(M(\alpha))=\mathbb{Z}_4\oplus \mathbb{Z}_4$.
Thus $M$ is a $\mathbb{Q}$-homology solid torus.
\end{proof}

\begin{lemma}\label{lem:244dash}
$G_P\cong H'(t/2,t,t)$ is impossible.
\end{lemma}

\begin{proof}
Each family of $t$ mutually parallel negative loops 
contains an $(i,i)$-edge for $i=1,2,\dots,t$.
Thus each vertex of $G_T$ is incident to two loops.
Then $G_T$ has $t/2$ components, each of which has an annulus support.
A jumping number argument easily rules out this.
\end{proof}

\bigskip

\textsc{Case (2)}. $p_0=t/2+1$.

\begin{lemma}\label{lem:finalcase}
$p_0=t/2+1$ is impossible.
\end{lemma}

\begin{figure}[tb]
\includegraphics*[scale=0.35]{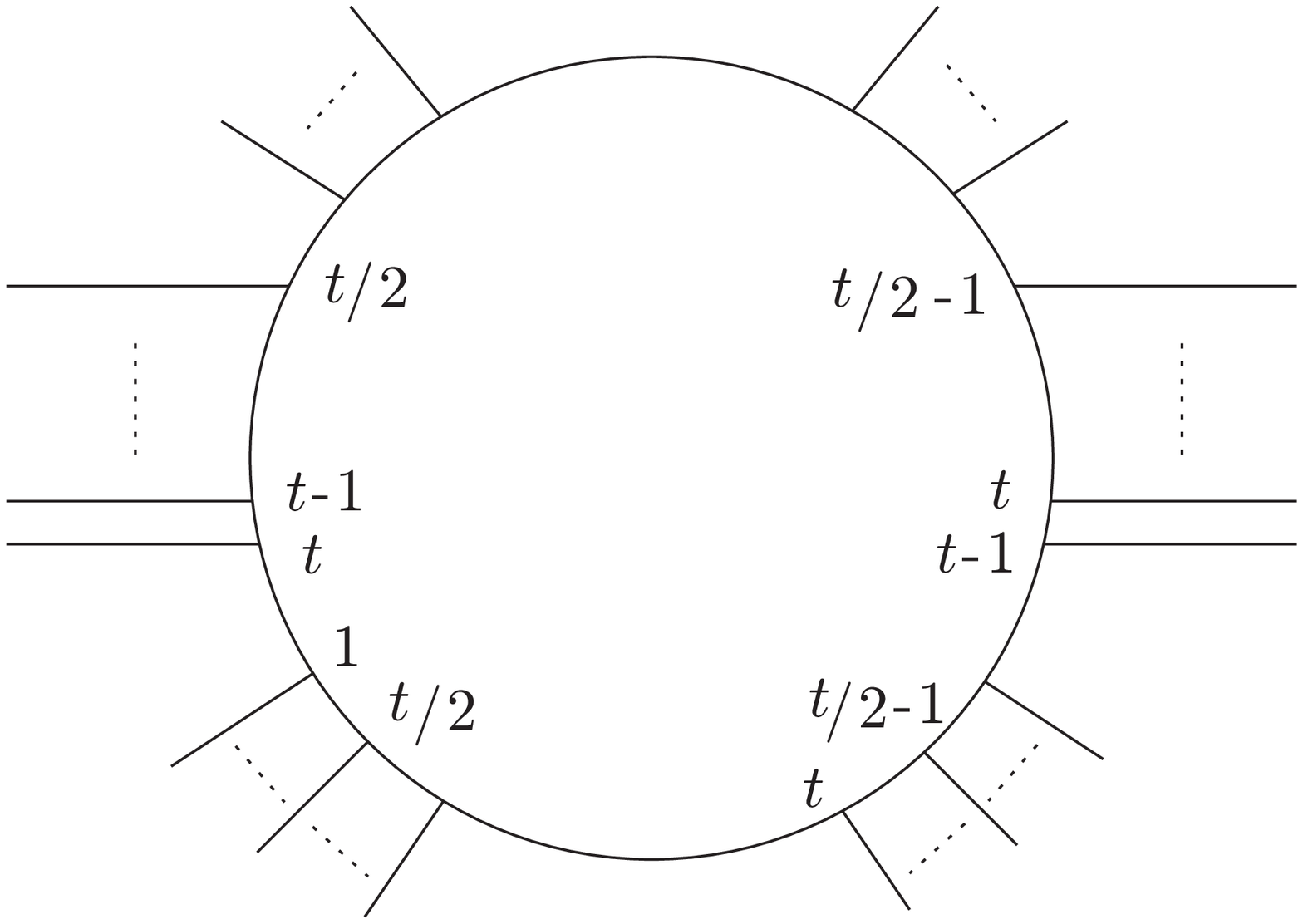}
\caption{}\label{fig:t/2+1}
\end{figure}

\begin{proof}
We see $G_P\cong H(t/2+1,t,t-1)$ or $H'(t/2+1,t,t-1)$.
For the former, $G_P$ contains a positive edge with the same label at its ends, a contradiction.
For the latter, we may assume that $G_P$ has labels as in Figure \ref{fig:t/2+1}.
Then $G_P$ contains a positive $(t/2-1,t/2)$-edge and a negative $(t/2-1,t/2)$-edge, a contradiction.
\end{proof}

\section{Proofs}\label{sec:proofs}

\begin{proof}[Proof of Theorem \ref{thm:main}]
Suppose that neither $M(\alpha)$ nor $M(\beta)$ contains a Klein bottle.
By Lemma \ref{lem:s2}, $s=2$.
If $t\le 2$, then Propositions \ref{pro:t1} and \ref{pro:t2} show that $M$ is a $\mathbb{Q}$-homology solid torus.
Proposition \ref{pro:tge3} rules out the case $t\ge 3$.

Suppose that $M(\alpha)$ or $M(\beta)$ contains a Klein bottle.
If both contain a Klein bottle, then $M=W(-4)$, which is a $\mathbb{Q}$-homology solid torus by Lemma \ref{lem:beta-noklein}.
For $t=1,2$, Proposition \ref{pro:t1kb}, Lemmas \ref{lem:t2p1} and \ref{lem:t2p2} give the conclusion.
If $t\ge 3$, then $p=1$ by Lemma \ref{lem:p1t3}.
By Lemma \ref{lem:positive-loop}, $G_P$ contains a positive loop, and
there are only three remaining cases after Lemma \ref{lem:finalcases}.
Then the remark after the proof of Lemma \ref{lem:finalcases}, Lemmas \ref{lem:244}, \ref{lem:244dash} and \ref{lem:finalcase}
lead us to the desired result.
We remark that if $p=1$ then $\partial N(\widehat{P})$ gives an essential torus in $M(\alpha)$ which meets the core of $V_\alpha$ in two points.
\end{proof}

The proof of Theorem \ref{thm:main} implicitly enables us to construct the collection of $3$-manifolds which includes
all hyperbolic $3$-manifolds with a single torus boundary $T_0$ such that there are two toroidal slopes on $T_0$ with distance $5$.
For example, consider the graph pair of Figure \ref{fig:t1s2-2}.
The manifold $X=N(S\cup T_0\cup T)$ has two $2$-spheres and $T_0$ as its boundary.
After capping the sphere components off with $3$-balls, we obtain a $3$-manifold $M$ with a single torus boundary $T_0$.
In this sense, we say that $M$ is uniquely determined from the graph pair.
But we do not know whether $M$ is hyperbolic or not.
The graph pairs of Figures \ref{fig:3331},
\ref{fig:kbt1}, \ref{fig:t2p1-0}, \ref{fig:22121}, and \ref{fig:244}
also determine the manifolds uniquely.
For the graph pairs of Figures \ref{fig:s2t2} and \ref{fig:t2p1}, the situation is different.
The manifold constructed as the above $X$ has the other torus component $T_1$ than $T_0$, after capping $2$-sphere components off with $3$-balls.
To obtain $M$ with a single torus boundary $T_0$, we need to perform Dehn filling on $T_1$.
Of course, there are infinitely many ways of Dehn fillings.
If we set $\mathcal{A}$ to be the collection of $3$-manifolds obtained from these $8$ graph pairs,
added $W(-4)$, where $W$ is the Whitehead link exterior,
then $\mathcal{A}$ includes all hyperbolic $3$-manifolds with a single torus boundary $T_0$
such that there are two toroidal slopes on $T_0$ with distance $5$.
Clearly, this collection consists of infinitely many manifolds.
For example, it contains the exteriors of Eudave-Mu\~{n}oz knots $k(2,-1,n,0)$ with $n\ne 1$.
(They correspond to the pair of Figure \ref{fig:s2t2} \cite{T2}.)

\begin{proof}[Proof of Corollary \ref{cor:main}]
If not, $\Delta(\alpha,\beta)=5$ \cite{Go2}.
By \cite{GT}, $\partial M$ is a single torus or two tori.
However, the latter case does not happen by \cite{L}.
For the former, Theorem \ref{thm:main} gives a contradiction.
\end{proof}

\begin{proof}[Proof of Corollary \ref{cor:hit}]
Let $\alpha$ and $\beta$ be two toroidal slopes on $T_0$.
Then $\Delta(\alpha,\beta)\le 8$ by \cite{Go2}.
Furthermore, if $\Delta(\alpha,\beta)=6$ or $8$, then
each of $M(\alpha)$ and $M(\beta)$ contains an essential torus
which meets the core of the attached solid torus in exactly two points.
When $\Delta(\alpha,\beta)=7$, one contains an essential torus meeting the core of the attached solid torus
in a single point, and the other contains one meeting the core of the attached solid torus in two points.

Suppose $\Delta(\alpha,\beta)=5$.
Then we showed that $\partial M$ is a single torus or two tori \cite{GT}.
By \cite{L}, the latter case happens only when $M$ is the Whitehead sister link exterior.
In this case, it is well known that each surgered manifold contains an essential torus meeting
the core of the attached solid torus just twice \cite{GW}.
Since the two components of the Whitehead sister link has non-zero linking number, $M$ cannot have
a properly embedded once-punctured torus.

If $\partial M$ is a single torus, Theorem \ref{thm:main} gives the conclusion.
\end{proof}


\bibliographystyle{amsplain}

\end{document}